\documentclass[11pt]{article}
\usepackage{amsfonts}
\usepackage{amssymb}
\newcommand{\be}{\begin{equation}}
\newcommand{\ee}{\end{equation}}
\newcommand{\ba}{\begin{array}}
\newcommand{\ea}{\end{array}}
\newcommand{\bea}{\begin{eqnarray}}
\newcommand{\eea}{\end{eqnarray}}
\newcommand{\bee}{\begin{eqnarray*}}
\newcommand{\eee}{\end{eqnarray*}}

\newtheorem{Thm}{Theorem}
\newtheorem{Lemma}{Lemma}
\newtheorem{Prop}{Proposition}

\newtheorem{Rk}{Remark}

\catcode`@=11
\renewcommand\appendix{\bigskip {\noindent\Large \bf Appendix}\par
  \setcounter{section}{0}%
  \setcounter{subsection}{0}%
  \renewcommand\thesection{\@Alph\c@section}}
\catcode`@=12

\def\thesection{\arabic{section}}

\def\C{{\bf C}}
\def\R{{\bf R}}

\def\lim{\mathop{\rm lim}}

\def\sup{\mathop{\rm sup}}

\def\e{\varepsilon}

\def\log{{\rm log}}

\def\scr{s_{c}}

\title{Blow up of the critical norm for some radial $L^2$ super critical non linear Schr\"odinger equations}\author{
Frank Merle\\
{\small  D\'epartement de Math\'ematiques}\\
{\small Universit\'e de Cergy-Pontoise}\\
{\small 2, avenue Adolphe Chauvin}\\
{\small 95302 Cergy-Pontoise Cedex France}\\
{\small frank.merle@math.u-cergy.fr}\\
 \ \\
Pierre Rapha\"el \\
{\small  Laboratoire de math\'ematiques}\\
{\small UMR 8628 du CNRS}\\
{\small Universit\'e Paris-Sud}\\ 
{\small 91405 Orsay Cedex France}\\
{\small pierre.raphael@math.u-psud.fr}}
\date{\today}

\begin{document}

\maketitle

\begin{abstract}

We consider the nonlinear Schr\"odinger equation $iu_t=-\Delta u-|u|^{p-1}u$ in dimension $N\geq 3$ in the $L^2$ super critical range $1+\frac{4}{N}< p<\frac{N+2}{N-2}$. The corresponding scaling invariant space is $\dot{H}^{s_c}$ with $0< s_c<1$ and this covers the physically relevant case $N=p=3$. The existence of finite time blow up solutions is known. Let $p_c=\frac{N}{2}(p-1)$ so that $\dot{H}^{s_c}\subset L^{p_c}$. Let $u(t)\in \dot{H}^{s_c}\cap \dot{H}^1$ be a radially symmetric blow up solution which blows up at $0<T<+\infty$, we prove that the scaling invariant $L^{p_c}$ norm also blows up with a lower bound: $$|u(t)|_{L^{p_c}}\geq |\log(T-t)|^{C_{N,p}} \ \ \mbox{as}  \ \ t\to T.$$
\end{abstract}


\section{Introduction}



\subsection{Setting of the problem}


We consider in this paper the nonlinear Schr\"odinger equation
\be
\label{nls}
(NLS) \ \  \left   \{ \begin{array}{ll}
         iu_t=-\Delta u-|u|^{p-1}u, \ \ (t,x)\in [0,T)\times \R^N\\
         u(0,x)=u_0(x), \ \ u_0:\R^N\to \C
         \end{array}
\right .
\ee
in dimension $N\geq 3$ with $$1<p<\frac{N+2}{N-2}.$$ From a result of Ginibre and Velo \cite{GV}, (\ref{nls}) is locally well-posed in $H^1=H^1(\R^N)$ and thus, for $u_0\in H^1$, there exists $0<T\leq +\infty$ and a unique solution $u(t)\in {\cal{C}}([0,T),H^1)$ to (\ref{nls}) and either $T=+\infty$, we say the solution is global, or $T<+\infty$ and then $\lim_{t\uparrow T}|\nabla u(t)|_{L^2}=+\infty$, we say the solution blows up in finite time.\\
(\ref{nls}) admits the following conservation laws in the energy space $H^1$: 
$$
\left   . \begin{array}{ll}
         L^2-\mbox{norm}: \ \ \int|u(t,x)|^2=\int|u_0(x)|^2;\\
         \mbox{Energy}:\ \ E(u(t,x))=\frac{1}{2}\int|\nabla u(t,x)|^2-\frac{1}{p+1}\int |u(t,x)|^{p+1}=E(u_0).
         \end{array}
\right .
$$
The scaling symmetry $\lambda^{\frac{2}{p-1}}u(\lambda^2t,\lambda x)$ leaves the homogeneous Sobolev space $\dot{H}^{s_c}$ invariant with 
\be
\label{defsc}
s_c=\frac{N}{2}-\frac{2}{p-1}.
\ee
It is classical from the conservation of the energy and the $L^2$ norm that for $s_c<0$, the equation is subcritical and all $H^1$ solutions are global and bounded in $H^1$. The smallest power for which blow up may occur is $p=1+\frac{4}{N}$ which corresponds to $s_c=0$ and is referred to as the $L^2$ critical case. The case $0<s_c< 1$ is the $L^2$ super critical and $H^1$ subcritical case.\\

We focus from now on onto the case $0\leq s_c<1$. The existence of finite time blow up solutions is a consequence of the virial identity, \cite{ZS}: let an initial condition $u_0\in \Sigma=H^1\cap\{xu\in L^2\}$ with $E(u_0)<0$, then the corresponding solution $u(t)$ to (\ref{nls}) satisfies $u(t)\in \Sigma$ with: 
\be
\label{virie}
\frac{d^2}{dt^2}\int|x|^2|u(t,x)|^2=4N(p-1)E(u_0)-\frac{16s_c}{N-2s_c}\int |\nabla u|^2\leq 16E(u_0)
\ee
and thus the positive quantity $\int|x|^2|u(t,x)|^2$ cannot exist for all times and $u$ blows up in finite time.\\

Recall now from Cazenave and Weissler \cite{CW} that given $u_0\in \dot{H}^{s_c}$, there exists a maximum time $T(u_0)>0$ and a unique maximal solution $u(t)\in{\cal{C}}([0,T(u_0)),\dot{H}^{s_c})$ to (\ref{nls}). Moreover, let the Strichartz pair $$\gamma=\frac{4(p+1)}{(p-1)(N-2s_c)}, \ \ \rho=\frac{N(p+1)}{N+s_c(p-1)},$$ then blowup is equivalent to 
\be
\label{cofohefho}
|u(t)|_{L^{\gamma}((0,T(u_0)),\dot{B}^{s_c}_{\rho,2})}=+\infty.
\ee 
More generally and following the same procedure, given $s_c<s\leq 1$ and $u_0\in \dot{H}^s$, there exists $T(s,u_0)>0$ and a unique maximal solution $u(t)\in{\cal{C}}([0,T(s,u_0)),\dot{H}^{s_c})$ to (\ref{nls}), and as the problem is now subcritical with respect to $\dot{H}^s$, $T(s,u_0)<+\infty$ iff $\lim_{t\to T(s,u_0)}|u(t)|_{\dot{H}^s}= +\infty$.\\

Let now $u_0\in \dot{H}^{s_c}\cap \dot{H}^1$, then there exists a maximum time $T>0$ and a unique maximal solution $u(t)\in{\cal{C}}([0,T),\dot{H}^{s_c})$ to (\ref{nls}). Indeed, the life times given by the local Cauchy theory in $\dot{H}^{s_c}$ and $\dot{H}^{s_c}\cap \dot{H}^1$ are the same from a standard argument -see Appendix A-. Moreover, if $u(t)$ blows up in finite time $0<T<+\infty$, then there holds the scaling lower bound:
\be
\label{scallinglowerbound}
\forall s_c<s\leq 1, \ \ |u(t)|_{\dot{H}^s}\geq \frac{C(N,p,s)}{(T-t)^{\frac{s-s_c}{2}}}.
\ee
Indeed, let $s_c<s\leq 1$, $t\in [0,T)$ and consider $v_t(\tau,x)=\lambda^{\frac{2}{p-1}}(t)u(t+\lambda^2(t)\tau,\lambda(t)x)$ with $\lambda^{s-s_c}(t)|u(t)|_{\dot{H}^s}=1$ so that $|v_t(0)|_{\dot{H}^s}=1$, then from the local Cauchy theory in $\dot{H}^s$ which is subcritical, there exists $\tau_0(s)>0$ such that $v$ is defined on $[0,\tau_0(s)]$ from which $t+\lambda^2(t)\tau_0(s)<T$, this is (\ref{scallinglowerbound}).\\

Let us remark that this argument does not apply for the critical $\dot{H}^{s_c}$ norm. On the basis of numerical simulations and formal arguments, it has been conjectured in the case $N=p=3$, see for example \cite{Co}, that at least for radial solutions, finite time blow implies:
\be
\label{efhnofi}
\lim_{t\to T}|u(t)|_{\dot{H}^{s_c}}=+\infty.
\ee
Note that this is in sharp contrast to the $L^2$ critical case $s_c=0$ where the $L^2$ norm is conserved and thus (\ref{efhnofi}) breaks down.\\

Note also that such kind of critical problems and behavior of the critical norms have been addressed in other settings, see for example Escauriaza, Seregin, Sverak \cite{ESS} for the 3D Navier-Stokes problem.


\subsection{A general strategy: reduction to a Liouville theorem}


Let us present a general and robust strategy to attack the proof of (\ref{efhnofi}) which is inspired form the works in  Martel, Merle \cite{MM2} and Merle, Rapha\"el \cite{MR3}. The idea is to first argue by contradiction and use compactness arguments to extract from a renormalized version of the solution an asymptotic object which generates a global in time nonpositive energy solution to (\ref{nls}). Here the arguments are quite general and could be extended to a wider class of solutions and problems. In a second step, one concludes using a Liouville type blow up result for nonpositive energy solutions.\\

More precisely, let $u_0\in \dot{H}^{s_c}\cap\dot{H}^1$ with radial symmetry and assume that the corresponding solution $u(t)$ to (\ref{nls}) blows up in finite time $0<T<\infty$ or equivalently from \cite{CW}: $$\lim_{t\to T}|\nabla u(t)|_{L^2}=+\infty.$$ Let a sequence $t_n\to T$ such that 
\be
\label{estuntn}
\lim_{t_n\to T}|\nabla u(t_n)|_{L^2}=+\infty
\ee 
and 
\be
\label{estdeuxtn}
\forall n\geq 1, \ \  |\nabla u(t_n)|_{L^2}=\max_{t\in [0,t_n]}|\nabla u(t)|_{L^2}.
\ee 
We now assume that:
\be
\label{hypcontrohsnorme}
\sup_{n\geq 1}|u(t_n)|_{\dot{H}^{s_c}}<+\infty
\ee
and look for a contradiction.\\
Let the sequence of rescaled initial data $$u_n(0,x)=\lambda_u(t_n)^{\frac{2}{p-1}}u(t_n,\lambda_u(t_n)x)$$ with 
\be
\label{calclulamndan}
\lambda_u(t)=\left(\frac{1}{|\nabla u(t)|_{L^2}}\right)^{\frac{1}{1-s_c}}\ \ \mbox{so that} \ \ |\nabla u_n(0)|_{L^2}=1.
\ee 
From the scaling invariance, (\ref{hypcontrohsnorme}) and the conservation of the energy, we have: 
\be
\label{estenergy}
|u_n(0)|_{\dot{H}^{s_c}}=|u(t_n)|_{\dot{H}^{s_c}}\leq C \ \ \mbox{and} \ \ E(u_n(0))=\lambda_u(t_n)^{2(1-s_c)}E(u_0)\to 0 
\ee
as $n\to +\infty$ and in particular
\be
\label{convfaibke}
u_n(0)\rightharpoonup v(0) \ \ \in\dot{H}^{s_c}\cap \dot{H}^1 \ \ \mbox{as} \ \ n\to +\infty
\ee
up to a subsequence. From the compact radial embedding $\dot{H}^{s_c}\cap \dot{H}^1\hookrightarrow L^{p+1}$, we have up to a subsequence $$u_n(0)\to v(0) \ \ \mbox{in} \ \ L^{p+1} \ \ \mbox{as} \ \ n\to +\infty$$ and thus $$E(v(0))\leq 0.$$ Observe now that the solution $u_n(\tau)$ to (\ref{nls}) with initial data $u_n(0)$ is explicitly 
\be
\label{hoziaaa}
u_n(\tau,x)=\lambda_u(t_n)^{\frac{2}{p-1}}u(t_n+\lambda_u(t_n)^2\tau,\lambda_u(t_n)x),
\ee 
and thus (\ref{estdeuxtn}) and (\ref{calclulamndan}) imply:
\be
\label{controlgefdaint}
\forall \tau \in (-\frac{t_n}{\lambda^2_u(t_n)},0], \ \ |\nabla u_n(\tau)|_{L^2}=\frac{|\nabla u(t_n+\lambda_u(t_n)^2\tau)|_{L^2}}{|\nabla u(t_n)|_{L^2}}\leq 1.
\ee 
Let now $v(t)$ be the solution to (\ref{nls}) with initial data $v(0)$ and $(-T_v,0]$ its maximum time interval existence on the left in time in $\dot{H}^{s_c}\cap \dot{H}^1$, then one may easily adapt the Lemma of stability of weak convergence in $H^1$, see Glangetas, Merle \cite{GM} and also Lemma 3 in \cite{MR3}, to conclude: $$\forall \tau<0, \ \ u_n(\tau)\rightharpoonup v(\tau) \ \ \mbox{in} \ \ \dot{H}^{s_c}\cap \dot{H}^1 \ \ \mbox{as} \ \ n\to +\infty$$ and thus from (\ref{controlgefdaint}): $$|\nabla v(\tau)|_{L^2}\leq 1 \ \ \mbox{and} \ \ -T_v=-\infty.$$ In other words, (\ref{hypcontrohsnorme}) implies the existence of a global in time radially symmetric nonpositive energy solution to (\ref{nls}) in $\dot{H}^{s_c}\cap \dot{H}^1$.\\ 

Our first result in this paper is that the existence of such an object may be ruled out in some cases from the following Liouville type result:

\begin{Thm}[Finite time blow up for non positive energy solutions in $\dot{H}^{s_c}\cap \dot{H}^1$]
\label{thm1}
Assume $$N\geq 3\ \  \mbox{and} \ \ 0<s_c<1.$$ Let $u_0\in \dot{H}^{\scr}\cap\dot{H}^1$ with radial symmetry and $$E(u_0)\leq 0,$$ then the corresponding solution $u(t)$ to (\ref{nls}) blows up in finite time $0<T<+\infty$.
\end{Thm}

{\bf Comments on Theorem \ref{thm1}}\\

{\it 1. On the assumption $u_0\in \dot{H}^{\scr}\cap\dot{H}^1$}: Let $u_0\in \Sigma=H^1\cap\{xu\in L^2\}$ with $E_0<0$, then finite time blow up follows from the virial identity (\ref{virie}). If $u_0\in H^1$ radial with $E_0<0$, a simple localization argument allows one to conclude also, see Ogawa, Tsutsumi \cite{OT}. Now if $u_0\in H^1$ radial with $E_0=0$, then finite time blow up also follows. The key here is first the conservation of the energy and a Gagliardo-Nirenberg inequality: $$|\nabla u(t)|_{L^2}^2=\frac{2}{p+1}|u|_{L^{p+1}}^{p+1}\leq C|\nabla u(t)|_{L^2}^{2+s_c(p-1)}|u(t)|_{L^2}^{(1-s_c)(p-1)}$$ which implies from {\it the conservation of the $L^2$ norm} the {\it uniform lower bound}: $$|\nabla u(t)|_{L^2}\geq \frac{C}{|u_0|_{L^2}^{\frac{1-s_c}{s_c}}}.$$ This together with a space localization of the virial identity (\ref{virie}) yields the claim.\\
Let us insist onto the fact that our need to work with low regularity $u_0\in \dot{H}^{\scr}\cap\dot{H}^1$ comes from the renormalization procedure before the extraction of the asymptotic object and a major difficulty is thus that we may no longer use the $L^2$ conservation law. Arguing by contradiction, we in fact need to rule out the possibility of a non linear self similar vanishing $|\nabla u(t)|_{L^2}\to 0$ as $t\to +\infty$ in the case when $E_0=0$. This difficulty already occurred in \cite{MR3}, \cite{MR4}. Our main tool is that for radial functions, we may replace the role of the $L^2$ norm by a suitable scaling invariant Morrey-Campanato norm, see the definition (\ref{defrho}), for which {\it uniform bounds in time are derived which somehow mimic the $L^2$ conservation law}. The key here is {\it a new kind of  monotonicity statement} based on a localized virial identity, see estimate (\ref{estrhouniform}) in Proposition \ref{lemmarhonorm}. Here the techniques are thus restricted to radial solutions.\\

{\it 2. On the sharpness of the result}: We expect the assumptions on the initial data to be sharp in the following sense. One may obtain exact self similar blow up solutions by looking for solutions of the form
\be
\label{seflsimilartype}
u(t,x)= \frac{1}{\lambda(t)^{\frac{2}{p-1}}}P\left(\frac{x}{\lambda(t)}\right)e^{i\log(T-t)} \ \ \mbox{with} \ \ \lambda(t)=\sqrt{2b(T-t)}
\ee 
for some  parameter $b>0$ and some stationary profile $P$ satisfying the non linear elliptic equation:
\be
\label{dhozei}
\Delta P -P + ib\left(\frac{2}{p-1}P+y\cdot\nabla P\right)+P|P|^{p-1}=0.
\ee
Rigorous existence results of finite energy radially symmetric solutions to (\ref{dhozei}) are known only for $p$ close to the $L^2$ critical value, see Rottsch\"afer and Kaper \cite{RK}. The obtained profiles are in $\dot{H}^1$ and have zero energy but always miss $\dot{H}^{s_c}$ due to a logarithmic growth at infinity. Such solutions, when they exist, thus provide explicit examples of zero energy solutions which blow up on the right in time but are global on the left.\\


\subsection{Lower bound on the blow up rate for the critical Sobolev norm}


We expect the above strategy to be quite robust. Moreover, it exhibits the fundamental baby problem which is the Liouville Theorem \ref{thm1}. It nevertheless has two weaknesses. First, it does not prove (\ref{efhnofi}) but only $$\limsup_{t\to T}|u(t)|_{\dot{H}^{s_c}}=+\infty,$$ a major difficulty being the possibility of oscillations in $|\nabla u(t)|_{L^2}$ which is a standard and difficult problem, see Martel and Merle \cite{MM3} for related problems for the generalized (KdV) equation, and \cite{MR1} for further discussions for the $L^2$ critical (NLS). Second, because the proof relies on an obstructive argument, it does not give any estimate on the rate of blow up of the critical Sobolev norm.\\

Our main claim in this paper is that both these difficulties may be overcome under the assumptions of Theorem \ref{thm1}. In other words, {\it for any sequence $t_n\to T$} and by pushing further the analysis of the proof of Theorem \ref{thm1} -giving an upper bound on the blow up time depending only on the repartition of the mass of the initial data in $\dot{H}^{s_c}$, see Remark \ref{rkk}-, we may {\it directly} estimate on the rescaled sequence $u_n(\tau)$ given by (\ref{hoziaaa}) the growth of the critical Sobolev norm and this provides us with a lower bound on its blow up rate.\\

In fact, we claim an even slightly stronger result:

\begin{Thm}[Lower bound for the critical $L^{p_c}$ norm]
\label{thmmain}
Assume $$N\geq 3\ \ \mbox{and} \ \ 0<s_c<1.$$ Let $$p_c=\frac{N}{2}(p-1) \ \ \mbox{so that} \ \ \dot{H}^{s_c}\subset L^{p_c}.$$  There exists a constant $\gamma=\gamma(N,p)>0$ such that the following holds true. Let $u_0\in \dot{H}^{\scr}\cap\dot{H}^1$ with radial symmetry and assume that the corresponding solution to (\ref{nls}) blows up in finite time $0<T<+\infty$, then 
\be
\label{lowerboundcritical}
|u(t)|_{L^{p_c}}\geq |\log(T-t)|^{\gamma}
\ee 
for $t$ close enough to $T$.
\end{Thm}

{\bf Comments on Theorem \ref{thmmain}}\\

{\it 1. Extension to $N=2$}: First observe from the Sobolev embedding $\dot{H}^{s_c}\subset L^{p_c}$ that (\ref{lowerboundcritical}) implies: 
\be
\label{norlecritiquelower}
|u(t)|_{\dot{H}^{s_c}}\geq C|\log(T-t)|^{\gamma}.
\ee
Moreover, the proof of both Theorem \ref{thm1} and Theorem \ref{thmmain} extends verbatim to the case $$N=2, \ \ 3<p<5.$$ The cases $p\geq 5$ in dimension $N=1,2$ remain open.\\

{\it 2. Sharpness of the result}: Let us remark that the physically relevant case $N=3=p=3$, $s_c=\frac{1}{2}$, $u_0\in H^1$, see \cite{SS}, is covered by Theorem \ref{thmmain}. The numerical simulations performed in this case, see \cite{SS} and references therein, strongly suggest the existence and the stability of a radially symmetric self similar blow up regime, that is with a blow up speed given by the scaling law $$|\nabla u(t)|_{L^2}\sim \frac{C(u_0)}{(T-t)^{\frac{1}{4}}}.$$ In this regime, we may derive an upper bound on the blow up rate of the critical $\dot{H}^{\frac{1}{2}}$ norm, and thus the $L^3$ norm, as follows. We write down Duhamel's formula and estimate using Strichartz estimates with the pair $(4,3)$, see \cite{S}, \cite{Cbook}:
$$|u(t)|_{\dot{H}^{\frac{1}{2}}} \leq   C|u(0)|_{\dot{H}^{\frac{1}{2}}}+C\left|\int_0^te^{i(t-\tau)\Delta}u(\tau)|u(\tau)|^2d\tau\right|_{L^{\infty}([0,t],\dot{B}^{\frac{1}{2}}_{2,2})}\leq  C(u_0)+|u|u|^2|_{L^{\frac{4}{3}}([0,t],\dot{B}^{\frac{1}{2}}_{\frac{3}{2},2})}.$$ We then have from Sobolev embeddings the standard nonlinear estimate in Besov spaces: $|u|u|^2|_{\dot{B}^{\frac{1}{2}}_{\frac{3}{2},2}}\leq C|u|_{L^6}^2|u|_{\dot{B}^{\frac{1}{2}}_{3,2}}\leq C|\nabla u|_{L^2}^3$, and thus $$|u(t)|_{L^3}\leq C|u(t)|_{\dot{H}^{\frac{1}{2}}}\leq C(u_0)+C\left(\int_0^t|\nabla u(\tau)|_{L^2}^4d\tau\right)^{\frac{3}{4}}\leq C|\log(T-t)|^{\frac{3}{4}}.$$ Combined with (\ref{lowerboundcritical}) and the explicit computation from the proof $\gamma(N=3,p=3)=\frac{1}{12}$, this gives $$|\log(T-t)|^{\frac{1}{12}}\leq |u(t)|_{L^3}\leq |\log(T-t)|^{\frac{3}{4}}.$$ We more generally expect that the lower bound (\ref{lowerboundcritical}) is sharp in the logarithmic scale, even though the constant $\gamma$ is most likely not sharp -and the value of the sharp constant is anyway unclear even at the formal level-.\\

{\it 2. On the lower bound (\ref{lowerboundcritical})}: From the proof, the lower bound (\ref{lowerboundcritical}) may be improved for solutions blowing up faster than the scaling estimate. More precisely, if one assumes $$|\nabla u(t)|_{L^2}\sim \frac{C}{(T-t)^{\alpha}} \ \ \mbox{with} \ \ \alpha>\frac{1-s_c}{2},$$ then the proof provides a lower bound $$|u(t)|_{L^{p_c}}\geq \frac{C}{(T-t)^{\beta}} \ \ \mbox{for some} \ \ \beta=\beta(\alpha,N,p)>0.$$ Note that the existence of such kind of excited solutions with respect to the scaling estimate is known only for $N=2$, $p=5$, see \cite{R2}, and indeed the critical norm blows up then at a polynomial rate. The existence of such solutions in the range of parameters of Theorem \ref{thmmain} is open.\\

The proof of Theorem \ref{thmmain} is a generalization of the one of Theorem \ref{thm1} and relies on fine properties of localization in space of the $L^2$ mass for the renormalized solution $u_n(\tau)$ given by (\ref{hoziaaa}), see Propositions \ref{lemmarhonorm} and \ref{lemmarhonormd}. In particular, let $t_n\to T$, we prove that the blow up of $u(t_n)$ implies a specific structure of $u_n(0)$ and more precisely a uniform lower bound on the scaling invariant weighted $L^2$ norm $$\int_{x\in {\cal{C}}_i}\frac{|u_n(0,x)|^2}{|x|^{2s_c}}dx\geq \frac{1}{|u(t_n)|^{\alpha}_{L^{p_c}}}$$ for some $\alpha=\alpha(N,p)>0$ and a well chosen family of disjoint annuli $\left({\cal{C}}_i\right)_{1\leq i\leq N(t)}$. Using H\"older inequality and summing over these annuli using an estimate $N(t)\geq |\log(T-t)|$ from the proof yields (\ref{lowerboundcritical}).\\

This paper is organized as follows. In section 2, we present some simple technical tools which we will need  along the proof. In section 3, we prove an abstract result for solutions to (\ref{nls}) which is the core of our analysis and is in fact a monotonicity type of statement, see Propositions \ref{lemmarhonorm} and \ref{lemmarhonormd}. In section 4, we first prove Theorem \ref{thm1} and then generalize the proof in order to get Theorem \ref{thmmain}.\\

{\bf Acknowledgments}: Both authors would like to thank the University of Chicago where part of this work was done. Pierre Rapha\"el would also like to thank the laboratoire MIP, Universit\'e Paul Sabatier, Toulouse, for its kind hospitality. In addition, we thank the referees for their careful reading of the paper.


\section{A radial interpolation estimate}


This section is devoted to the proof of a radial interpolation estimate needed for the proof of Theorems \ref{thmmain} and \ref{thm1}.\\

Let us recall the standard Gagliardo-Nirenberg inequality 
\be
\label{galiehozinc}
\int|u|^{p+1}\leq C|u|_{L^{p_c}}^{p-1}\left(\int |\nabla u|^2\right).
\ee  
For radially symmetric distributions, we may sharpen this inequality away from the origin using the $\rho$ semi-norm 
\be
\label{defrho}
\rho(u,R)=\sup_{R'\geq R}\frac{1}{(R')^{2s_c}}\int_{R'\leq |x|\leq 2R'}|u|^2.
\ee and the fact that for $N\geq 3$ and $s_c<1$, $$p<\frac{N+2}{N-2}\leq 5,$$ and thus the nonlinearity is $L^2$ subcritical away from zero. 

\begin{Lemma}[Radial Gagliardo-Nirenberg inequality]
\label{lemmacleinterpol}
(i) There exists a universal constant $C>0$ such that for all $u\in L^{p_c}$, 
\be
\label{estunfluxrho}
\forall R>0, \ \ \frac{1}{R^{2s_c}}\int_{|y|\leq R}|u|^2\leq C|u|^2_{L^{p_c}}
\ee
 and 
\be
\label{estecayzero}
\lim_{R\to +\infty} \frac{1}{R^{2s_c}}\int_{|x|\leq R}|u|^2\to 0 \ \ \mbox{as} \ \ R\to +\infty.
\ee
(ii) For all $\eta>0$, there exists a constant $C_{\eta}>0$ such that for all $u\in \dot{H}^{s_c}\cap \dot{H}^1$ with radial symmetry, for all $R>0$, 
\be
\label{ineginterpol}
 \int_{|x|\geq R}|u|^{p+1}\leq \eta|\nabla u|_{L^2(|x|\geq R)}^2+\frac{C_{\eta}}{R^{2(1-s_c)}}\left[(\rho(u,R))^{\frac{p+3}{5-p}}+(\rho(u,R))^{\frac{p+1}{2}}\right].
\ee
\end{Lemma}

{\bf Proof of Lemma \ref{lemmacleinterpol}}\\

(i) From H\"older: $$\frac{1}{R^{2s_c}}\int_{|y|\leq R}|u|^2\leq C|u|_{L^{p_c}(|y|\leq R)}^2$$ and (\ref{estunfluxrho}) and (\ref{estecayzero}) follow -the last one by splitting the integral in two-.\\ 
(ii) We split the proof of (\ref{ineginterpol}) in two steps:\\

{\bf step 1} Localized interpolation estimate on a ring.\\

Let a smooth radially symmetric $u$ and pick an $\eta>0$. Consider for $D>0$ the annulus ${\cal{C}}=\{D\leq |x|\leq 2D\}$, we claim: 
\be
\label{estkeyzero}
\int_{{\cal{C}}}|u|^{p+1}\leq \eta|\nabla u|^{2}_{L^2({\cal{C}})}+\frac{C_{\eta}}{D^{2(1-s_c)}}\left[(\rho(u,D))^{\frac{p+3}{5-p}}+(\rho(u,D))^{\frac{p+1}{2}}\right].
\ee
Indeed, let $x_0\in {\cal{C}}$ such that $|u(x_0)|=|u|_{L^{\infty}({\cal{C}})}$. Then either there exists $y\in {\cal{C}}$ such that 
\be
\label{conduy}
|u(y)|\leq \frac{1}{2} |u(x_0)|
\ee 
in which case $$|u|^2_{L^{\infty}({\cal{C}})}\leq C||u(x_0)|^2-|u(y)|^2|\leq \int_{D}^{2D}|u||u'|\leq \frac{C}{D^{N-1}}|\nabla u|_{L^2({\cal{C}})}|u|_{L^2({\cal{C}})}$$ and thus $$|u|_{L^{\infty}({\cal{C}})}\leq \frac{C}{D^{\frac{N-1}{2}}}|\nabla u|^{\frac{1}{2}}_{L^2({\cal{C}})}|u|^{\frac{1}{2}}_{L^2({\cal{C}})}.$$ This implies 
\bee
\int_{{\cal{C}}}|u|^{p+1} & \leq & |u|^{p-1}_{L^{\infty}({\cal{C}})}\int_{{\cal{C}}}|u|^2\leq \frac{C}{D^{\frac{(N-1)(p-1)}{2}}}|\nabla u|^{\frac{p-1}{2}}_{L^2({\cal{C}})}|u|_{L^2({\cal{C}})}^{\frac{p+3}{2}}\\
& \leq & C|\nabla u|^{\frac{p-1}{2}}_{L^2({\cal{C}})}\left(\frac{1}{D^{2s_c}}\int_{D\leq |x|\leq 2D}|u|^2\right)^{\frac{p+3}{4}}\frac{1}{D^{\frac{(N-1)(p-1)}{2}-\frac{2(p+3)s_c}{4}}}\\
& \leq & C|\nabla u|^{\frac{p-1}{2}}_{L^2({\cal{C}})}\left[\rho(u,D)\right]^{\frac{p+3}{4}}\frac{1}{D^{\frac{(5-p)(1-s_c)}{2}}}
\eee
where we used from direct computation
$$\frac{(N-1)(p-1)}{2}-\frac{2(p+3)s_c}{4}=\frac{(5-p)(1-s_c)}{2}.$$ (\ref{estkeyzero}) now follows from $\frac{p-1}{2}<2$ and H\"older.\\
Now if (\ref{conduy}) never holds on ${\cal{C}}$, then $$|u|_{L^{\infty}({\cal{C}})}\leq \frac{C}{D^{\frac{N}{2}}}|u|_{L^2({\cal{C}})}$$ from which 
\bee
& & \int_{{\cal{C}}}|u|^{p+1}  \leq   D^N |u|_{L^{\infty}({\cal{C}})}^{p+1}  \leq  \frac{C}{D^{\frac{N(p+1)}{2}-N}}|u|_{L^2({\cal{C}})}^{p+1}\\
& \leq & \frac{1}{D^{\frac{N(p+1)}{2}-N-s_c(p+1)}}\left(\frac{1}{D^{2s_c}}\int_{D\leq |x|\leq 2D}|u|^2\right)^{\frac{p+1}{2}}\leq  \frac{C}{D^{2(1-s_c)}}\left[\rho(u,D)\right]^{\frac{p+1}{2}}
\eee
where we used from explicit computation $$\frac{N(p+1)}{2}-N-s_c(p+1)=2(1-s_c),$$ and (\ref{estkeyzero}) holds true again.\\

{\bf step 2} Conclusion.\\

Given $R>0$, we write $$\int_{|x|\geq R}|u|^{p+1}=\Sigma_{j=0}^{+\infty}\int_{2^jR\leq |x|\leq 2^{j+1}R}|u|^{p+1}$$ and apply (\ref{estkeyzero}) with $D=R2^j$. Note from the monotonicity of $\rho(u,R)$ that $$\forall j\geq 0, \ \ \rho(u,2^jR)\leq \rho(u,R)$$ from which 
\bee
\int_{|x|\geq R}|u|^{p+1} & \leq & \eta\Sigma_{j=0}^{+\infty}|\nabla u|^2_{L^2(R2^j\leq |x|\leq 2^{j+1}R)}+C_{\eta}\left[(\rho(u,R))^{\frac{p+3}{5-p}}+(\rho(u,R))^{\frac{p+1}{2}}\right]\Sigma_{j=0}^{+\infty}\frac{1}{(R2^j)^{2(1-s_c)}}\\
& \leq & \eta|\nabla u|_{L^2(|x|\geq R)}^2+\frac{C_{\eta}}{R^{2(1-s_c)}}\left[(\rho(u,R))^{\frac{p+3}{5-p}}+(\rho(u,R))^{\frac{p+1}{2}}\right].
\eee
This concludes the proof of Lemma \ref{lemmacleinterpol}.


\section{The main Propositions}


In this section, we prove the main Propositions at the heart of the proof of Theorem \ref{thmmain}.\\

Let us rewrite the Gagliardo-Nirenberg inequality (\ref{galiehozinc}) as follows: there exists a universal constant $C_{GN}=C_{GN}(N,p)>0$ such that:
\be
\label{estggalibdeux}
\forall u\in \dot{H}^{s_c}\cap \dot{H}^1, \ \ E(u)\geq \frac{1}{2}\int |\nabla u|^2\left[1-\left(\frac{|u|_{L^{p_c}}}{C_{GN}}\right)^{p-1}\right].
\ee
We first claim a uniform control of the scaling invariant $\rho$ norm and a dispersive control of radially symmetric solutions to (\ref{nls}) in a parabolic region in time. Indeed we prove under suitable assumptions on the initial data -$\dot{H}^{s_c}$ control and a priori bound on non positive part of the Hamiltonian- that the $\dot{H}^1$ norm must decay in a self similar fashion in average.

\begin{Prop}[Uniform control of the $\rho$ norm and dispersion]
\label{lemmarhonorm}
Let $N\geq 3$ and $0<s_c<1$. There exist universal constants $C_1,\alpha_1,\alpha_2>0$ depending on $N,p$ only such that the following holds true. Let $v(\tau)\in {\cal{C}}([0,\tau_*], \dot{H}^{s_c}\cap \dot{H}^1)$ be a radially symmetric solution to (\ref{nls}) and assume:
\be
\label{estuimatestzero}
\tau_*^{1-s_c}\max(E(v_0),0)<1,
\ee
and 
\be
\label{estsizemzero}
M_0=\frac{4|v(0)|_{L^{p_c}}}{C_{GN}}\geq 2
\ee
where $C_{GN}$ is the universal constant in the Gagliardo-Nirenberg inequality (\ref{estggalibdeux}). Then there holds the uniform control of the $\rho$ norm:
\be
\label{estrhouniform}
\rho(v(\tau_*),M_0^{\alpha_1}\sqrt{\tau_*})\leq C_1M_0^2,
\ee
and the global dispersive estimate:
\be
\label{dispeestimategradient}
\int_0^{\tau_*}(\tau_*-\tau)|\nabla v(\tau)|_{L^2}^2d\tau\leq M_0^{\alpha_2}\tau_*^{1+s_c}.
\ee
\end{Prop}

{\bf Proof of Proposition \ref{lemmarhonorm}}\\

{\bf step 1} Localized radial virial estimate.\\

Let a smooth radially symmetric cut-off function function $\psi(x)=\frac{|x|^2}{2}$ for $|x|\leq 2$ and $\psi(x)=0$ for $|x|\geq 3$ with 
\be
\label{contorlderiveeut}
|\nabla \psi|^2\leq C\psi,
\ee 
We claim that there exist constants $C'_1(N,p), C'_2(N,p)>0$ such that $\forall R>0$, $\forall \tau\in [0,\tau_*]$, 
\bea
\label{virialcontrolocalized}
 \nonumber & &  C'_1\int |\nabla v|^2 +\frac{d}{d\tau}Im\left(\int\nabla \psi_R\cdot\nabla v\overline{v}\right)\\
& \leq  & N(p-1)E(v_0)+C_2'\left[\int_{|x|\geq R}|v|^{p+1}+\frac{1}{R^2}\int_{2R\leq |x|\leq 3R}|v|^2\right].
\eea
Proof of (\ref{virialcontrolocalized}): Let $\chi$ be a smooth radially symmetric compactly supported cut-off function. We recall the following standard localized virial identities which up to standard regularization arguments are obtained by integration by parts on (\ref{nls}):
\be
\label{estfluxldeux}
\frac{1}{2}\frac{d}{d\tau}\int \chi|v|^2=Im\left(\int\nabla \chi\cdot\nabla v\overline{v}\right),
\ee
\be
\label{estirialalgebraic}
\frac{1}{2}\frac{d}{d\tau}Im\left(\int\nabla \chi\cdot\nabla v\overline{v}\right)=\int\chi''|\nabla v|^2-\frac{1}{4}\int \Delta^2\chi |v|^2-\left(\frac{1}{2}-\frac{1}{p+1}\right)\int\Delta\chi|v|^{p+1}.
\ee
Note that we used here that $v$ has radial symmetry. We apply (\ref{estirialalgebraic}) with $\chi(x)=\psi_R(x)=R^2\psi(\frac{x}{R})$ and get:
\bee
\nonumber & & \frac{1}{2}\frac{d}{d\tau}Im\left(\int\nabla \psi_R\cdot\nabla v\overline{v}\right)\\
& = & \int\psi''(\frac{x}{R})|\nabla v|^2-\frac{1}{4R^2}\int \Delta^2\psi(\frac{x}{R}) |v|^2-\left(\frac{1}{2}-\frac{1}{p+1}\right)\int\Delta\psi(\frac{x}{R})|v|^{p+1}\\
\nonumber & \leq & \int |\nabla v|^2-N\left(\frac{1}{2}-\frac{1}{p+1}\right)\int|v|^{p+1}+C\left[\frac{1}{R^2}\int_{2R\leq |x|\leq 3R}|v|^2+\int_{|x|\geq R}|v|^{p+1}\right].
\eee
From the conservation of the energy $\int|v|^{p+1}=\frac{p+1}{2}\int |\nabla v|^2-(p+1)E(v_0)$ from which 
$$
\int |\nabla v|^2-N\left(\frac{1}{2}-\frac{1}{p+1}\right)\int|v|^{p+1}=\frac{N(p-1)}{2}E(v_0)-\frac{2s_c}{N-2s_c}\int|\nabla v|^2,
$$
and (\ref{virialcontrolocalized}) follows.\\

{\bf step 2} A priori control of the $\rho$ norm on parabolic space time intervals.\\

First recall from (\ref{estunfluxrho}), (\ref{estsizemzero}) that there exists a universal constant $C_2>0$ such that 
\be
\label{estinitialzero}
\forall R>0, \ \ \frac{1}{R^{2s_c}}\int_{|x|\leq 3R}|v(0)|^2\leq C_2M_0^2.
\ee 
We now claim the following a priori control: there exists a universal constant $C>0$ such that for all $A>0$ and $\tau_0\in [0,\tau^*]$, let 
\be
\label{defmtildeifnity}
M^2_{\infty}(A,\tau_0)=\max_{\tau\in [0,\tau_0]}\rho(v(\tau), A\sqrt{\tau}), \ \ R=A\sqrt{\tau_0},
\ee 
then: 
\bea
\label{newestimateone}
\nonumber & & C_1'\int_0^{\tau_0}(\tau_0-\tau)|\nabla v(\tau)|_{L^2}^2d\tau \leq  C\tau_0^{1+s_c}\left[M_0^2A^{2(1+s_c)}+ \frac{[M_{\infty}(A,\tau_0)]^{\frac{2(p+3)}{5-p}}+[M_{\infty}(A,\tau_0)]^{2}}{A^{2(1-s_c)}}\right]\\
& + & 2\tau_0\left[Im\int\nabla \psi_R\cdot\nabla v(0)\overline{v(0)}+\frac{N(p-1)E(v_0)}{2}\tau_0\right],
\eea
\bea
\label{newestimateoneflxuldeux}
\nonumber & & \frac{1}{R^{2s_c}}\int_{R\leq |x|\leq 2R}|v(\tau_0)|^2 \leq  8C_2M_0^2+\frac{C}{A^4}\left[[M_{\infty}(A,\tau_0)]^{\frac{2(p+3)}{5-p}}+[M_{\infty}(A,\tau_0)]^2\right]\\
& + & \frac{4}{\tau_0^{s_c}A^{2(1+s_c)}}\left[Im\int\nabla \psi_R\cdot\nabla v(0)\overline{v(0)}+\frac{N(p-1)E(v_0)}{2}\tau_0\right].
\eea
Proof of (\ref{newestimateone}) and (\ref{newestimateoneflxuldeux}): Consider the estimate (\ref{virialcontrolocalized}) for the localization $$R=A\sqrt{\tau_0}$$ and estimate the terms of the right hand side. From the definition (\ref{defrho}) of $\rho$ and the definition (\ref{defmtildeifnity}) of $M_{\infty}(A,\tau_0)$, we have: $$\forall \tau\in [0,\tau_0], \ \ \rho(v(\tau),R)=\rho(v(\tau),A\sqrt{\tau_0})\leq \rho(v(\tau),A\sqrt{\tau})\leq M^2_{\infty}(A,\tau_0).$$ We thus estimate:
\be
\label{estldeuxun}
\frac{1}{R^2}\int_{2R\leq   |x|\leq 3R}|v(\tau)|^2   \leq   \frac{C}{R^{2(1-s_c)}}\rho(v(\tau),R)\leq \frac{CM^2_{\infty}(A,\tau_0)}{R^{2(1-s_c)}}
\ee
The non linear term in (\ref{virialcontrolocalized}) is estimated from the refined Gagliardo-Nirenberg estimate (\ref{ineginterpol}) provided $\eta>0$ has been chosen small enough: $\forall \tau\in[0,\tau_0]$,
\bea
\label{estldeuxtrois}
\nonumber C'_2\int_{|x|\geq R}|v|^{p+1} & \leq &  \frac{C'_1}{2}\int |\nabla v(\tau)|^2+\frac{C}{R^{2(1-s_c)}}\left[(\rho(v(\tau),R))^{\frac{p+3}{5-p}}+(\rho(v(\tau),R))^{\frac{p+1}{2}}\right]\\
& \leq & \frac{C'_1}{2}\int |\nabla v(\tau)|^2+C\frac{[M_{\infty}(A,\tau_0)]^{\frac{2(p+3)}{5-p}}+[M_{\infty}(A,\tau_0)]^{2}}{R^{2(1-s_c)}}
\eea
where we used $\frac{p+3}{5-p}>\frac{p+1}{2}>1$. We now inject (\ref{estldeuxun}) and (\ref{estldeuxtrois}) into (\ref{virialcontrolocalized}) and integrate in time to get: $\forall \tau\in[0,\tau_0]$, 
\bee
& & \frac{C_1'}{2}\int_0^{\tau}|\nabla v(\sigma)|^2_{L^2}d\sigma +Im\left(\int \nabla \psi_R\cdot\nabla v(\tau)\overline{v(\tau)}\right)\\
&  \leq  & Im\left(\int\nabla \psi_R\cdot\nabla v(0)\overline{v(0)}\right)+ N(p-1)E(v_0)\tau+C\frac{[M_{\infty}(A,\tau_0)]^{\frac{2(p+3)}{5-p}}+[M_{\infty}(A,\tau_0)]^{2}}{R^{2(1-s_c)}}\tau.
\eee 
We integrate once more in time between $0$ and $\tau_0$ from (\ref{estfluxldeux}) and get: 
\bea
\label{estindknhfuie}
\nonumber  & & \int \psi_{R}|v(\tau_0)|^2+C_1'\int_0^{\tau_0}(\tau_0-\tau)|\nabla v(\tau)|^2_{L^2}d\tau\\
\nonumber & \leq & \int \psi_{R}|v(0)|^2+  2\tau_0\left[Im\left(\int\nabla \psi_R\cdot\nabla v(0)\overline{v(0)}\right)+\frac{N(p-1)E(v_0)}{2}\tau_0\right]\\
\nonumber &+& C\frac{\tau_0^2}{R^{2(1-s_c)}}\left[[M_{\infty}(A,\tau_0)]^{\frac{2(p+3)}{5-p}}+[M_{\infty}(A,\tau_0)]^{2}\right]\\
\nonumber &\leq &  4C_2M_0^2R^{2(1+s_c)}+ 2\tau_0\left[Im\int\nabla \psi_R\cdot\nabla v(0)\overline{v(0)}+\frac{N(p-1)E(v_0)}{2}\tau_0\right ]\\
& + & C\tau_0^{1+s_c}\left[\frac{[M_{\infty}(A,\tau_0)]^{\frac{2(p+3)}{5-p}}+[M_{\infty}(A,\tau_0)]^{2}}{A^{2(1-s_c)}}\right]
\eea
where we used in the last step (\ref{estinitialzero}) and the definition of $R$ (\ref{defmtildeifnity}). This with (\ref{estuimatestzero}) and (\ref{defmtildeifnity}) implies (\ref{newestimateone}).\\
We now divide (\ref{estindknhfuie}) by $R^{2(1+s_c)}$ and observe from the definition of $\psi$ that:
$$\frac{1}{R^{2(1+s_c)}}\int \psi_{R}|v(\tau_0)|^2\geq \frac{1}{2R^{2s_c}}\int_{R\leq |x|\leq 2R}|v(\tau_0)|^2,$$ and (\ref{newestimateoneflxuldeux}) follows from the definition of $R$ (\ref{defmtildeifnity}) again.\\

{\bf step 3} Self similar decay of the gradient.\\

We now claim that one may derive the global dispersive estimate (\ref{dispeestimategradient}) as a consequence of (\ref{newestimateone}) and (\ref{newestimateoneflxuldeux}).\\

Indeed, let $\e>0$ be a small enough constant depending only on $(N,p)$ to be chosen later. Let 
\be
\label{defgnot}
G_{\e}=M_0^{\frac{1}{\e}}, \ \ A_{\e}=\left(\frac{\e G_{\e}}{M_0^2}\right)^{\frac{1}{2(1+s_c)}}.
\ee
Recall from (\ref{estunfluxrho}) that there exists a universal constant $C>0$ such that $$\forall R>0, \ \ \forall u\in L^{p_c}, \ \ \rho(u,R)\leq C|u|_{L^{p_c}}^2.$$ From the regularity of the flow $v\in{\cal{C}}([0,\tau^*],\dot{H}^{s_c}\cap\dot{H}^1)$, the Sobolev embedding $\dot{H}^{s_c}\subset L^{p_c}$ and the definition (\ref{estsizemzero}) of $M_0$, we may consider the largest time $\tau_1\in [0,\tau_*]$ such that 
\be
\label{estdisperiveboosteddeux}
\forall \tau_0\in[0,\tau_1], \ \ \int_0^{\tau_0}(\tau_0-\tau)|\nabla v(\tau)|_{L^2}^2d\tau\leq G_{\e}\tau_0^{1+s_c}
\ee
and
\be
\label{initiationdepartdeux}
M^2_{\infty}(A_{\e},\tau_1)=\max_{\tau\in [0,\tau_1]}\rho(v(\tau), A_{\e}\sqrt{\tau})\leq 2\frac{M_0^2}{\e}.
\ee
We claim that:
\be
\label{estdisperiveboosteddeuxboot}
\forall \tau_0\in[0,\tau_1], \ \ \int_0^{\tau_0}(\tau_0-\tau)|\nabla v(\tau)|_{L^2}^2d\tau\leq \frac{G_{\e}}{2}(1+\tau_0)^{1+s_c},
\ee
and 
\be
\label{initiationdepartdeuxboot}
M^2_{\infty}(A_{\e},\tau_1)=\max_{\tau\in [0,\tau_1]}\rho(v(\tau), A_{\e}\sqrt{\tau})\leq \frac{M_0^2}{\e}
\ee
provided $\e>0$ has been chosen small enough, and (\ref{estrhouniform}) and (\ref{dispeestimategradient}) follow.\\
Proof of (\ref{estdisperiveboosteddeuxboot}) and (\ref{initiationdepartdeuxboot}): Let $\tau_0\in [0,\tau_1]$. We rewrite (\ref{newestimateone}) with $A=A_{\e}$ and $R=A_{\e}\sqrt{\tau_0}$ and estimate the terms in the right hand side. Observe from (\ref{estsizemzero}), (\ref{defgnot}) and (\ref{initiationdepartdeux}) that: 
\bea
 \label{cnfeohfroro}
\nonumber & &\frac{[M_{\infty}(A_{\e},\tau_0)]^{\frac{2(p+3)}{5-p}}+[M_{\infty}(A_{\e},\tau_0)]^{2}}{A_{\e}^{2(1-s_c)}} \leq \left(\frac{M_0}{\e}\right)^C\frac{1}{G_{\e}^{\frac{1-s_c}{1+s_c}}}\\
& \leq & \frac{1}{\e^CM_0^{{\frac{1-s_c}{1+s_c}\frac{1}{\e}}-C}}\leq \frac{1}{\e^C 2^{\frac{1}{C\e}}} \leq \frac{1}{10}
\eea
 for $\e$ small enough. Injecting this into (\ref{newestimateone}) and using (\ref{defgnot}) yields:
\bea
\label{newestimateonebisbis}
 \nonumber & & C_1'\int_0^{\tau_0}(\tau_0-\tau)|\nabla v(\tau)|_{L^2}^2d\tau \leq  C\tau_0^{1+s_c}\left[M_0^2A_{\e}^{2(1+s_c)}+ \frac{[M_{\infty}(A_{\e},\tau_0)]^{\frac{2(p+3)}{5-p}}+[M_{\infty}(A_{\e},\tau_0)]^{2}}{A_{\e}^{2(1-s_c)}}\right]\\
\nonumber & + & 2\tau_0\left[Im\int\nabla \psi_R\cdot\nabla v(0)\overline{v(0)}+\frac{N(p-1)E(v_0)}{2}\tau_0\right]\\
\nonumber & \leq & C\tau_0^{1+s_c}\left[\e G_{\e}+\frac{1}{10}\right]+ 2\tau_0\left[Im\int\nabla \psi_R\cdot\nabla v(0)\overline{v(0)}+\frac{N(p-1)E(v_0)}{2}\tau_0\right]\\
& \leq & G_{\e}\tau_0^{1+s_c}\left\{\frac{1}{10}+\frac{2}{G_{\e}\tau_0^{s_c}}\left[Im\int\nabla \psi_R\cdot\nabla v(0)\overline{v(0)}+\frac{N(p-1)E(v_0)}{2}\tau_0\right]\right\},
\eea
provided $\e>0$ has been chosen small enough. We now claim the following key estimate which controls the growth of the momentum term in the above right hand side: $\forall \tau_0\in[0,\tau_1]$, $\forall A\geq A_{\e}$, let $R=A\sqrt{\tau_0}$, then:
\be
\label{keycontrolmomentum}
Im\int\nabla \psi_R\cdot\nabla v(0)\overline{v(0)}+\frac{N(p-1)E(v_0)}{2}\tau_0\leq C\frac{M_0^2A^{2(1+s_c)}}{\e^{\frac{1}{1+s_c}}}\tau_0^{s_c}.
\ee
Let us assume (\ref{keycontrolmomentum}) and conclude the proof. We inject  (\ref{keycontrolmomentum}) at $A=A_{\e}$ with (\ref{defgnot}) into (\ref{newestimateonebisbis}) and get: 
\bee
& & \int_0^{\tau_0}(\tau_0-\tau)|\nabla v(\tau)|_{L^2}^2d\tau \leq G_{\e}\tau_0^{1+s_c}\left[\frac{1}{10}+ C\frac{M_0^2A_{\e}^{2(1+s_c)}}{G_{\e}\e^{\frac{1}{1+s_c}}}\right]\\
& = & G_{\e}\tau_0^{1+s_c}\left[\frac{1}{10}+ C\frac{\e}{\e^{\frac{1}{1+s_c}}}\right]= G_{\e}\tau_0^{1+s_c}\left[\frac{1}{10}+ C\e^{\frac{s_c}{1+s_c}}\right]\leq \frac{G_{\e}}{2}\tau_0^{1+s_c}
\eee
provided $\e>0$ has been chosen small enough, and (\ref{estdisperiveboosteddeuxboot}) follows. Moreover, let $A\geq A_{\e}$ and consider (\ref{newestimateoneflxuldeux}) for $R=A\sqrt{\tau_0}$. First observe that $\rho$ is a non increasing function of $R$ so that from (\ref{initiationdepartdeux}): $$M_{\infty}(A,\tau_0)\leq M_{\infty}(A_{\e},\tau_0)\leq M_{\infty}(A_{\e},\tau_1)\leq \frac{M_0^2}{\e}.$$ We inject this together with (\ref{keycontrolmomentum}) into (\ref{newestimateoneflxuldeux}) and argue like for the proof of (\ref{cnfeohfroro}) to derive:
\bee
\nonumber & & \frac{1}{R^{2s_c}}\int_{R\leq |x|\leq 2R}|v(\tau_0)|^2  \leq  8C_2M_0^2+\frac{C}{A^4}\left[[M_{\infty}(A,\tau_0)]^{\frac{2(p+3)}{5-p}}+[M_{\infty}(A,\tau_0)]^2\right]\\
& + & \frac{4}{A^{2(1+s_c)}\tau_0^{s_c}}\left[Im\int\nabla \psi_R\cdot\nabla v(0)\overline{v(0)}+\frac{N(p-1)E(v_0)}{2}\tau_0\right]\\
&  \leq  &  M_0^2\left[4C_2+\frac{1}{\e^CM_0^{\frac{1}{C\e}}}+\frac{C}{\e^{\frac{1}{1+s_c}}}\right]\leq \frac{M_0^2}{\e}
\eee
for $\e>0$ small enough, and (\ref{initiationdepartdeuxboot}) is proved.\\

{\bf step 4} Proof of the momentum estimate (\ref{keycontrolmomentum}).\\

We now turn to the proof of (\ref{keycontrolmomentum}).\\
Let us first make an observation and assume for simplicity that $E(v_0)\leq 0$. From $R=A\sqrt{\tau}_0$,  (\ref{keycontrolmomentum}) is implied by: $$\left|Im\int\nabla \psi_R\cdot\nabla v(0)\overline{v(0)}\right|\leq C\frac{M_0^2A^{2}}{\e^{\frac{1}{1+s_c}}}R^{2s_c}.$$ Remarkably enough, an estimate of this type is straightforward for $\frac{1}{2}\leq s_c<1$ using the following interpolation estimate: $$\forall v_0\in \dot{H}^{s_c}, \ \ \forall R>0, \ \ \left|Im\int\nabla \psi_R\cdot\nabla v(0)\overline{v(0)}\right|\leq C|v_0|_{\dot{H}^{s_c}}^2R^{2s_c}.$$ This would provide a short cut to the proof of (\ref{norlecritiquelower}) for $\frac{1}{2}\leq s_c<1$. For $0<s_c\leq \frac{1}{2}$, such an estimate however is certainly not true for any $v_0$ and the best one could hope for is: $$ \forall v_0\in \dot{H}^{s_c}, \ \ \forall R>0, \ \ \left|Im\int\nabla \psi_R\cdot\nabla v(0)\overline{v(0)}\right|\leq CR|\nabla v_0|_{L^2}^{\frac{1-2s_c}{1-s_c}}|v_0|_{\dot{H}^{s_c}}^{\frac{1}{1-s_c}}$$ which is not enough to conclude.\\ 
We now come back to the proof of (\ref{keycontrolmomentum}) and claim that it indeed holds true as a consequence of the backwards integration of the flow of (\ref{nls}) on parabolic regions in time in the dispersive regime (\ref{estdisperiveboosteddeux}). Let $\tau_0\in [0,\tau_1]$, $A\geq A_{\e}$ and $R=A\sqrt{\tau_0}$.\\
We first claim that there exists a universal constant $K(N,p)>0$ and 
\be
\label{cbdbonol}
\tilde{\tau}_0\in [\frac{\e^{\frac{1}{1+s_c}}}{4}\tau_0,\frac{\e^{\frac{1}{1+s_c}}}{2}\tau_0] \ \ \mbox{with} \ \ |\nabla v(\tilde{\tau}_0)|^2_{L^2}\leq \frac{KG_{\e}}{\tilde{\tau_0}^{1-s_c}}.
\ee 
Proof of (\ref{cbdbonol}): If not, let $\tilde{\tau}=\e^{\frac{1}{1+s_c}}\tau_0$, then $$\int_{\frac{\tilde{\tau}}{4}}^{\frac{\tilde{\tau}}{2}}|\nabla v(\sigma)|_{L^2}^2d\sigma\geq KG_{\e}\int_{\frac{\tilde{\tau}}{4}}^{\frac{\tilde{\tau}}{2}}\frac{d\sigma}{\sigma^{1-s_c}}\geq CKG_{\e}\tilde{\tau}^{s_c}.$$ On the other hand, we have $\tilde{\tau}=\e^{\frac{1}{1+s_c}}\tau_0\leq \tau_0\leq \tau_1$ and thus from (\ref{estdisperiveboosteddeux}): $$G_{\e}\tilde{\tau}^{1+s_c}\geq \int_0^{\tilde{\tau}}(\tilde{\tau}-\sigma)|\nabla v(\sigma)|_{L^2}^2d\sigma\geq \frac{\tilde{\tau}}{2}\int_{\frac{\tilde{\tau}}{4}}^{\frac{\tilde{\tau}}{2}}|\nabla v(\sigma)|_{L^2}^2d\sigma\geq CKG_{\e}\tilde{\tau}^{1+s_c}$$ and a contradiction follows for $K>0$ large enough. This concludes the proof of (\ref{cbdbonol}).\\
Let now 
\be
\label{cofehoh}
R=A\sqrt{\tau_0}=A_1\sqrt{\tilde{\tau}_0} \ \ \mbox{and thus} \ \ \frac{\e^{\frac{1}{2(1+s_c)}}}{16}\leq \frac{A}{A_1}\leq \e^{\frac{1}{2(1+s_c)}}.
\ee
We claim: 
\be
\label{cankzws}
\left|Im\int\nabla \psi_R\cdot\nabla v(\tilde{\tau}_0)\overline{v(\tilde{\tau}_0)}\right|\leq \frac{CM_0^2A^{2(1+s_c)}}{\e^{\frac{1}{1+s_c}}}\tau_0^{s_c}.
\ee
Proof of (\ref{cankzws}): We first claim that:
\be
\label{nvkfngvepriog}
\frac{1}{R^{2(1+s_c)}}\int \psi_R|v(\tilde{\tau}_0)|^2\leq C\left[M_0^2+\frac{G_{\e}}{A_1^{2(1+s_c)}}\right].
\ee
Indeed, consider (\ref{estfluxldeux}) with $\chi(x)=\psi_R(x)=R^2\psi(\frac{x}{R})$, then from Cauchy-Schwarz and (\ref{contorlderiveeut}):
$$\left|\frac{d}{d\tau}\int \psi_R |v|^2\right|=2\left|Im\left(\int\nabla \psi_R\cdot\nabla v\overline{v}\right)\right|\leq C|\nabla v|_{L^2}\left(\int \psi_R|v|^2\right)^{\frac{1}{2}}.$$ Integrating this differential inequation between $0$ and $\tilde{\tau}_0$ and using (\ref{estinitialzero}) yields: 
\bea
\label{ciquoouporut}
\nonumber \int\psi_R|v(\tilde{\tau}_0)|^2   & \leq &  C\left[\int\psi_R|v(0)|^2+\left(\int_0^{\tilde{\tau}_0}|\nabla v(\sigma)|_{L^2}d\sigma\right)^2\right]\\
& \leq & C\left[R^{2(1+s_c)}M_0^2+\tilde{\tau}_0\int_0^{\tilde{\tau}_0}|\nabla v(\sigma)|^2_{L^2}d\sigma\right].
\eea
We now observe from (\ref{cbdbonol}) that $2\tilde{\tau}_0\leq \e^{\frac{1}{1+s_c}}\tau_0\leq \tau_0\leq\tau_1$  and thus (\ref{estdisperiveboosteddeux}) ensures:
\be
\label{ionguefgiitildetuao}
\int_0^{\tilde{\tau}_0}|\nabla v(\sigma)|^2_{L^2}d\sigma \leq \frac{1}{\tilde{\tau}_0}\int_0^{2\tilde{\tau}_0}(2\tilde{\tau}_0-\sigma)|\nabla v(\sigma)|^2_{L^2}d\sigma\leq CG_{\e}\tilde{\tau_0}^{s_c}.
\ee Injecting this into (\ref{ciquoouporut}) and using (\ref{cofehoh}) yields: $$\int\psi_R|u(\tilde{\tau}_0)|^2\leq C\left[R^{2(1+s_c)}M_0^2+G_{\e}\tilde{\tau}_0^{1+s_c}\right]=CR^{2(1+s_c)}\left[M_0^2+\frac{G_{\e}}{A_1^{2(1+s_c)}}\right]$$ and concludes the proof of (\ref{nvkfngvepriog}).\\
We now estimate from (\ref{cbdbonol}) and (\ref{nvkfngvepriog}):
\bee
 & & \left|Im\int\nabla \psi_R\cdot\nabla v(\tilde{\tau}_0)\overline{v(\tilde{\tau}_0)}\right|\leq R^{1+s_c}|\nabla v(\tilde{\tau}_0)|_{L^2}\left(\frac{1}{R^{2(1+s_c)}}\int \psi_R|v(\tilde{\tau}_0)|^2\right)^{\frac{1}{2}}\\
& \leq & CR^{1+s_c}\frac{G_{\e}^{\frac{1}{2}}}{\tilde{\tau_0}^{\frac{1-s_c}{2}}}\left[M_0^2+\frac{G_{\e}}{A_1^{2(1+s_c)}}\right]^{\frac{1}{2}}\leq CA^{1+s_c}\tau_0^{\frac{1+s_c}{2}}\frac{G_{\e}^{\frac{1}{2}}}{(\tau_0\frac{A^2}{A^2_1})^{\frac{1-s_c}{2}}}\left[M_0+\frac{G_{\e}^{\frac{1}{2}}}{A_1^{(1+s_c)}}\right]\\
& = & CM_0^2A^{2(1+s_c)}\tau_0^{s_c}\left[\left(\frac{G_{\e}}{A^{2(1+s_c)}M_0^2}\right)^{\frac{1}{2}}\left(\frac{A_1}{A}\right)^{1-s_c}+\frac{G_{\e}}{A^{2(1+s_c)}M_0^2}\left(\frac{A}{A_1}\right)^{2s_c}\right]\\
& \leq & CM_0^2A^{2(1+s_c)}\tau_0^{s_c}\left[\left(\frac{G_{\e}}{A_{\e}^{2(1+s_c)}M_0^2}\right)^{\frac{1}{2}}\left(\frac{A_1}{A}\right)^{1-s_c}+\frac{G_{\e}}{A_{\e}^{2(1+s_c)}M_0^2}\left(\frac{A}{A_1}\right)^{2s_c}\right]
\eee
from $A\geq A_{\e}$. We now inject (\ref{defgnot}) and (\ref{cofehoh}) to get:
\bee
\nonumber  & & \left|Im\int\nabla \psi_R\cdot\nabla v(\tilde{\tau}_0)\overline{v(\tilde{\tau}_0)}\right|\\
& \leq & CM_0^2A^{2(1+s_c)}\tau_0^{s_c}\left[\frac{1}{\e^{\frac{1}{2}}\e^{\frac{1-s_c}{2(1+s_c)}}}+\frac{\e^{\frac{2s_c}{2(1+s_c)}}}{\e}\right]= \frac{2CM_0^2A^{2(1+s_c)}}{\e^{\frac{1}{1+s_c}}}\tau_0^{s_c}.
\eee
This concludes the proof of (\ref{cankzws}).\\
We now integrate the control (\ref{cankzws}) backwards at $\tau=0$. To wit, consider (\ref{estirialalgebraic}) with $\chi=\psi_R$ and estimate the terms in the right hand side. From (\ref{defgnot}), (\ref{initiationdepartdeux}):
\bee
\left|\int \Delta^2\psi_R |v|^2\right| &\leq  & \frac{C}{R^2}\int_{2R\leq |x|\leq 3R}|v(\tau)|^2\leq \frac{CM^2_{\infty}(A_{\e},\tau_1)}{R^{2(1-s_c)}}\\
& \leq  & \frac{1}{\tilde{\tau}_0^{1-s_c}}\frac{G_{\e}^{C\e}}{A_{\e}^{2(1-s_c)}}\leq \frac{1}{\tilde{\tau}_0^{1-s_c}}
\eee
for $\e>0$ small enough. We then use the conservation of the energy to derive from (\ref{estirialalgebraic}) the crude estimate: 
$$
\left|\frac{d}{d\tau}Im\left(\int\nabla \psi_R\cdot\nabla v\overline{v}\right)\right|\leq C\left[\int|\nabla v|^2+|E(v_0)|+\frac{1}{\tilde{\tau}_0^{1-s_c}}\right]
$$
We integrate this in time from $0$ to $\tilde{\tau}_0$ and get:
$$
\left|Im\int\nabla \psi_R\cdot\nabla v(0)\overline{v(0)}\right|\leq \left|Im\int\nabla \psi_R\cdot\nabla v(\tilde{\tau}_0)\overline{v(\tilde{\tau}_0)}\right|+C\left[\int_0^{\tilde{\tau}_0}|\nabla v(\sigma)|_{L^2}^2d\sigma+|E(v_0)|\tilde{\tau}_0+\tilde{\tau}^{s_c}_0\right].
$$
This implies in particular using (\ref{cankzws}) and (\ref{ionguefgiitildetuao}):
\bea
\label{cbjzoodo}
\nonumber & & Im\int\nabla \psi_R\cdot\nabla v(0)\overline{v(0)}+\frac{N(p-1)E(v_0)}{2}\tau_0  \leq   \frac{CM_0^2A^{2(1+s_c)}}{\e^{\frac{1}{1+s_c}}}\tau_0^{s_c}+CG_{\e}\tilde{\tau}_0^{s_c}\\
& + & \frac{N(p-1)E(v_0)}{2}\tau_0+C|E(v_0)|\tilde{\tau}_0.
\eea
It remains to estimate each term in the above right hand side. For the second term, we have from (\ref{defgnot}) and (\ref{cbdbonol}): 
\be
\label{cbhoeoefoh}
G_{\e}\tilde{\tau}_0^{s_c}\leq CA_{\e}^{2(1+s_c)}M_0^2\frac{\e^{\frac{s_c}{1+s_c}}}{\e}\tau_0^{s_c}\leq \frac{CM_0^2A^{2(1+s_c)}}{\e^{\frac{1}{1+s_c}}}\tau_0^{s_c}.
\ee
Eventually, the third term in (\ref{cbjzoodo}) is estimated using (\ref{estuimatestzero}) and (\ref{cbdbonol}):
\bee
\frac{N(p-1)E(v_0)}{2}\tau_0+C\tilde{\tau}_0|E(v_0)|& \leq & \left[\frac{N(p-1)E(v_0)}{2}+C\e^{\frac{1}{1+s_c}}|E_0|\right]\tau_0\\
& \leq & C\max(E(v_0),0)\tau_0\leq C\tau_0^{s_c}
\eee 
for $\e>0$ small enough and where we used (\ref{defgnot}) in the last step. Injecting this with (\ref{cbhoeoefoh}) into (\ref{cbjzoodo}) concludes the proof of (\ref{keycontrolmomentum}).\\

This concludes the proof of Proposition \ref{lemmarhonorm}.\\

Under the hypothesis of Proposition \ref{lemmarhonorm}, we assume an additional energetic constraint on the solution and prove from the control of the $\rho$ norm (\ref{estrhouniform}) and the dispersive estimate (\ref{dispeestimategradient}) that it implies a non trivial repartition of the $L^2$ mass of the initial data. 

\begin{Prop}[Lower bound on a weighted local $L^2$ norm of $v(0)$]
\label{lemmarhonormd}
Let $N\geq 3$ and $0<s_c<1$. There exist universal constants $\alpha_3,c_3>0$ depending on $N,p$ only such that the following holds true. Let $v(\tau)\in {\cal{C}}([0,\tau_*], \dot{H}^{s_c}\cap \dot{H}^1)$ be a radially symmetric solution to (\ref{nls}) such that (\ref{estuimatestzero}) and (\ref{estsizemzero}) of Proposition \ref{lemmarhonorm} hold. Let
\be
\label{deflambdrescalied}
\lambda_v(\tau)=\left(\frac{1}{|\nabla v(\tau)|_{L^2}}\right)^{\frac{1}{1-s_c}}.
\ee 
Let
\be
\label{esttautildestar}
\tau_0\in [0,\frac{\tau_*}{2}]
\ee
and assume moreover the energetic constraint:
\be
\label{hypsupplemenatire}
\lambda^{2(1-s_c)}_v(\tau_0)E(v_0)\leq \frac{1}{4}.
\ee
Let
\be
\label{defdzerofun}
F_*=\frac{\sqrt{\tau_0}}{\lambda_v(\tau_0)}
\ee
and 
\be
\label{defdzero}
D_*=M_0^{\alpha_3}\max\left[1, F_*\right],
\ee
then 
\be
\label{estsclaingldeux}
\frac{1}{\lambda^{2s_c}_v(\tau_0)}\int_{|x|\leq D_*\lambda_v(\tau_0)}|v(0)|^2\geq c_3.
\ee
\end{Prop}

{\bf Proof of Proposition \ref{lemmarhonormd}}\\

{\bf step 1} Energetic constraint and lower bound on $v(\tau_0)$.\\

We claim that the energetic constraint (\ref{hypsupplemenatire}) together with the uniform control of the $\rho$ norm (\ref{estrhouniform}) imply (\ref{estsclaingldeux}) for $v(\tau_0)$. More precisely, we claim that there exist universal constants $C_3,c_3>0$ such that:
\be
\label{estsclaingldeuxstar}
\frac{1}{\lambda^{2s_c}_v(\tau_0)}\int_{|x|\leq A_*\lambda_v(\tau_0)}|v(\tau_0)|^2\geq c_3 
\ee
with
\be
\label{defcdeepslionb}
A_*=C_3\max[M_0^{\alpha_1}F_*, M_0^{\frac{p+3}{(1-s_c)(5-p)}}].
\ee
Proof of (\ref{estsclaingldeuxstar}): Consider a renormalization of $v(\tau_0)$: 
\be
\label{dehw}
w(x)=\lambda^{\frac{2}{p-1}}_v(\tau_0)v(\tau_0,\lambda_v(\tau_0)x),
\ee 
then from (\ref{deflambdrescalied}), (\ref{hypsupplemenatire}) and the conservation of the energy: 
\be
\label{ehxshazdio}
|\nabla w|_{L^2}=1 \ \ \mbox{and} \ \ E(w)=\lambda^{2(1-s_c)}_v(\tau_0)E(v(\tau_0))=\lambda^{2(1-s_c)}_v(\tau_0)E(v_0)\leq \frac{1}{4}
\ee 
and thus 
\be
\label{estlepdfonamental}
\int |w|^{p+1}=(p+1)\left(\frac{1}{2}\int |\nabla w|^2-E(w)\right)\geq \frac{p+1}{4}.
\ee
Pick now $\e>0$ small enough and let 
\be
\label{deaespilon}
A_{\e}=C_{\e}\max[M_0^{\alpha_1}F_*, M_0^{\frac{p+3}{(1-s_c)(5-p)}}]
\ee 
for $C_{\e}$ large enough to be chosen. First observe from (\ref{estrhouniform}), (\ref{defdzerofun}), (\ref{dehw}), (\ref{deaespilon}) and the scaling invariance of the $\rho$ norm that: 
\bea
\label{estrhorecales}
\nonumber \rho(w,A_{\e}) & \leq &  \rho(w,M_0^{\alpha_1}F_*)= \rho(v(\tau_0), M_0^{\alpha_1}\lambda_v(\tau_0)F_*)\\
& = &   \rho(v(\tau_0), M_0^{\alpha_1}\sqrt{\tau_0})\leq C_1M_0^2.
\eea
Thus, from (\ref{ineginterpol}) and (\ref{deaespilon}):
\bea
\label{estkeyun}
\nonumber \int_{|x|\geq A_{\e}}|w|^{p+1} & \leq & \e|\nabla w|_{L^2}^2+\frac{\hat{C}_{\e}}{A_{\e}^{2(1-s_c)}}\left[(\rho(w,A_{\e}))^{\frac{p+3}{5-p}}+(\rho(w,A_{\e}))^{\frac{p+1}{2}}\right]\\
& \leq & \e+\frac{\hat{C}_{\e}}{A_{\e}^{2(1-s_c)}}M_0^{\frac{2(p+3)}{5-p}}\leq 2\e
\eea
for $C_{\e}>1$ large enough. Injecting this into (\ref{estlepdfonamental}) yields for $\e>0$ small enough:$$ \int_{|x|\leq A_{\e}}|w|^{p+1}\geq \frac{p+1}{8}.$$ $\e>0$ being now fixed, recall the Gagliardo-Nirenberg inequality $$\int|w|^{p+1}\leq C|w|_{L^2}^{(p-1)(1-s_c)}|\nabla w|_{L^2}^{2+(p-1)s_c}$$ which we may localize to get: 
\bee
\frac{p+1}{8} & \leq &  \int_{|y|\leq A_{\e}}|w|^{p+1}\leq C|w|_{L^2(|y|\leq 2A_{\e})}^{(p-1)(1-s_c)}|\nabla w|_{L^2}^{2+(p-1)s_c}+C|w|_{L^2(|y|\leq 2A_{\e})}^{p+1}\\
& \leq & C|w|_{L^2(|y|\leq 2A_{\e})}^{(p-1)(1-s_c)}+C|w|_{L^2(|y|\leq 2A_{\e})}^{p+1}
\eee
and thus $$\int_{|y|\leq 2A_{\e}} |w|^2\geq c_3>0$$ for some constant $c_3=c_3(N,p)>0$. But from (\ref{dehw}), this is $$\frac{1}{\lambda^{2s_c}_v(\tau_0)}\int_{|x|\leq A_{\e}\lambda_v(\tau_0)}|v(\tau_0)|^2\geq c_3,$$ and this concludes the proof of (\ref{estsclaingldeuxstar}).\\

{\bf step 2} Backwards integration of the $L^2$ fluxes.\\

We now integrate backwards the $L^2$ fluxes from $\tau_0$ to $0$. The key here is that the global dispersive estimate (\ref{dispeestimategradient}) implies that such kind of local $L^2$ quantities are {\it almost conserved in time} in the parabolic region. We claim: $\forall \e>0$, there exists $\tilde{C}_{\e}>0$ such that $\forall D\geq D_{\e}$ with 
\be
\label{defcdeepslion}
D_{\e}=\tilde{C}_{\e}F_*\max\left[M_0^{\alpha_1}, M_0^{\frac{\alpha_2}{2(1+s_c)}}\right],
\ee
let 
\be
\label{choixdern}
\tilde{R}=\tilde{R}(D,\tau_0)=D\lambda_v(\tau_0),
\ee 
and $\chi_{\tilde{R}}(r)=\chi(\frac{r}{\tilde{R}})$ for some smooth radially symmetric cut-off function $\chi(r)=1$ for $r\leq 1$, $\chi(r)=0$ for $r\geq 2$, then:
\be
\label{estfluxldeuxun}
\left|\frac{1}{\lambda^{2s_c}_v(\tau_0)}\int \chi_{\tilde{R}}|v(\tau_0)|^2-\frac{1}{\lambda^{2s_c}_v(\tau_0)}\int\chi_{\tilde{R}}|v(0)|^2\right|<\e.
\ee
(\ref{estsclaingldeuxstar}) and (\ref{estfluxldeuxun}) now imply (\ref{estsclaingldeux}).\\
Proof of (\ref{estfluxldeuxun}): Pick $\e>0$. We compute the $L^2$ fluxes from (\ref{estfluxldeux}) with $\chi_{\tilde{R}}$:
\bea
\label{cnenpo}
\nonumber \left|\frac{d}{d\tau}\int \chi_{\tilde{R}}|v|^2\right| &= & 2\left|Im\left(\int\nabla \chi_{\tilde{R}}\cdot\nabla v\overline{v}\right)\right|\leq \frac{C}{\tilde{R}}|\nabla v(\tau)|_{L^2}\left(\int_{\tilde{R}\leq |x|\leq 2\tilde{R}}|v(\tau)|^2\right)^{\frac{1}{2}}\\
& \leq & \frac{C}{\tilde{R}^{1-s_c}}|\nabla v(\tau)|_{L^2}\left(\rho(v(\tau),\tilde{R})\right)^{\frac{1}{2}}.
\eea
Now observe from (\ref{defdzerofun}), (\ref{defcdeepslion}) and (\ref{choixdern}) that:
$$\forall \tau\in [0,\tau_0], \ \ \tilde{R}=D\lambda_v(\tau_0)\geq D_{\e}\frac{\lambda_v(\tau_0)}{\sqrt{\tau_0}}\sqrt{\tau_0}\geq \frac{D_{\e}}{F_*}\sqrt{\tau}\geq M_0^{\alpha_1}\sqrt{\tau},$$ and thus (\ref{estrhouniform}) and the monotonicity of $\rho$ ensure:
\be
\label{estkeyrehomonotonie}
\forall \tau\in[0,\tau_0], \ \ \rho(v(\tau),\tilde{R})\leq \rho(v(\tau),M_0^{\alpha_1}\sqrt{\tau})<CM_0^2.
\ee
Injecting this into (\ref{cnenpo}) yields: $$\forall \tau\in [0,\tau_0], \ \ \left|\frac{d}{d\tau}\int \chi_{\tilde{R}}|v(\tau)|^2\right|\leq \frac{CM_0}{\tilde{R}^{1-s_c}}|\nabla v(\tau)|_{L^2}.$$ We integrate this between $0$ and $\tau_0$, divide by $\tilde{R}^{2s_c}$ and use (\ref{dispeestimategradient}) and (\ref{esttautildestar}) to get:
\bee
 & & \left|\frac{1}{\lambda^{2s_c}_v(\tau_0)}\int \chi_{\tilde{R}}|v(\tau_0)|^2-\frac{1}{\lambda^{2s_c}_v(\tau_0)}\int\chi_{\tilde{R}}|v(0)|^2\right|\\
& \leq & \frac{C}{\tilde{R}^{1+s_c}}\int_0^{\tau_0}|\nabla v(\tau)|_{L^2}d\tau\leq \frac{C}{\tilde{R}^{1+s_c}}\left(\tau_0\int_0^{\tau_0}|\nabla v(\tau)|_{L^2}^2d\tau\right)^{\frac{1}{2}}\\
& \leq & \frac{C}{\tilde{R}^{1+s_c}}\left(\int_0^{2\tau_0}(2\tau_0-\tau)|\nabla v(\tau)|_{L^2}^2d\tau\right)^{\frac{1}{2}}\leq CM_0^{\frac{\alpha_2}{2}}\left(\frac{\tau_0}{\tilde{R}^2}\right)^{\frac{1+s_c}{2}}\\
& = & CM_0^{\frac{\alpha_2}{2}}\left(\frac{F_*}{D}\right)^{1+s_c}\leq CM_0^{\frac{\alpha_2}{2}}\left(\frac{F_*}{D_{\e}}\right)^{1+s_c}.
\eee
(\ref{estfluxldeuxun}) now follows from (\ref{defcdeepslion}) for $\tilde{C}_{\e}>1$ in (\ref{defcdeepslion}) large enough. \\

This concludes the proof of Proposition \ref{lemmarhonormd}.


\section{Proof of the main theorems}


This section is devoted to the proof of the main theorems. We start by proving the Liouville Theorem \ref{thm1}, and then refine the corresponding analysis to derive the proof of Theorem \ref{thmmain} and of the lower bound (\ref{lowerboundcritical}). 


\subsection{Proof of the Liouville Theorem \ref{thm1}}


Let $u_0\in \dot{H}^{\scr}\cap\dot{H}^1$ with radial symmetry and $E(u_0)\leq 0$, we argue by contradiction and assume that $u(\tau)$ is globally defined on the right $u\in {\cal{C}}([0,+\infty), \dot{H}^{s_c}\cap\dot{H}^1).$ First observe from $E_0\leq 0$ and (\ref{estggalibdeux}) that $u_0$ satisfies (\ref{estuimatestzero}) and (\ref{estsizemzero}), and thus (\ref{dispeestimategradient}) of Proposition \ref{lemmarhonorm} implies: 
$$\forall \tau_*\geq 0, \ \ \int_0^{\tau_*}(\tau_*-\tau)|\nabla u(\tau)|_{L^2}^2d\tau\leq C(u_0)\tau_*^{1+s_c}.
$$
This implies the existence of a sequence $\tau_n\to +\infty$ such that $$\forall n\geq 0, \ \ |\nabla u(\tau_n)|_{L^2}\leq \frac{C(u_0)}{\tau_n^{\frac{1-s_c}{2}}}$$ or equivalently from (\ref{deflambdrescalied}): 
\be
\label{nezolnmnc}
\lambda_u(\tau_n)=\left(\frac{1}{|\nabla u(\tau_n)|_{L^2}}\right)^{\frac{1}{1-s_c}}\geq C(u_0)\sqrt{\tau_n}.\ee
Apply now Proposition \ref{lemmarhonormd} at time $\tau_0=\tau_n$. Observe that (\ref{hypsupplemenatire}) is fulfilled from $E(u_0)\leq 0$. Moreover, from (\ref{nezolnmnc}), $F_n=F_*$ given by (\ref{defdzerofun}) satisfies:
$$\forall n\geq 0, \ \ F_n=\frac{\sqrt{\tau_n}}{\lambda_u(\tau_n)}\leq C(u_0).$$ From (\ref{defdzero}), we conclude that there exists $D_*$ {\it independent of $n$} such that (\ref{estsclaingldeux}) holds ie: 
\be
\label{chioefhapzeopze}
\frac{1}{\lambda^{2s_c}_u(\tau_n)}\int_{|x|\leq D_*\lambda_u(\tau_n)}|u(0)|^2\geq c_3>0.
\ee 
But from (\ref{nezolnmnc}), $\lambda_u(\tau_n)\to +\infty$ as $n\to +\infty$ and thus $u_0\in \dot{H}^{s_c}\subset L^{p_c}$ and (\ref{estecayzero}) contradict (\ref{chioefhapzeopze}). This concludes the proof of Theorem \ref{thm1}.

\begin{Rk}
\label{rkk} Note that in fact, the same argument allows one to derive in the setting of Theorem \ref{thm1} an explicit upper estimate on the blow up time of the solution depending only on the distribution function of the mass of the initial data $g(R)=\frac{1}{R^{2s_c}}\int_{|x|\leq R}|u(0)|^2$ and its size in $L^{p_c}$.
\end{Rk}


\subsection{Proof of Theorem \ref{thmmain}}


Let $u_0\in \dot{H}^{\scr}\cap\dot{H}^1$ with radial symmetry and assume that the corresponding solution $u(t)$ to (\ref{nls}) blows up in finite $0<T<+\infty$. Pick $t$ close enough to $T$. Let 
\be
\label{deflambdaut}
\lambda_u(t)=\left(\frac{1}{|\nabla u(t)|_{L^2}}\right)^{\frac{1}{1-s_c}},
\ee
 then the scaling lower bound (\ref{scallinglowerbound}) implies:
\be
\label{chohoeih}
\lambda_u(t)\leq C(N,p)\sqrt{T-t}.
\ee
Let the renormalization of $u(t)$: 
\be
\label{defvrescale}
v^{(t)}(\tau,x)=\lambda^{\frac{2}{p-1}}_u(t)\overline{u}(t-\lambda^2_u(t)\tau,\lambda_u(t)x).
\ee
To clarify the notations, we omit the dependence of $v$ on $t$. Note that $v(\tau)\in {\cal{C}}([0,\tau_{max}], \dot{H}^{\scr}\cap\dot{H}^1)$ is a solution to (\ref{nls}) with $\tau_{max}=\frac{t}{\lambda_u^2(t)}$. Let 
\be
\label{defnt}
N(t)=-\log\lambda_u(t) \ \ \mbox{so that} \ \ \tau_{max}=te^{2N(t)}.
\ee
Our aim is to prove:
\be
\label{tobeprovedcrucial}
|v(0)|_{L^{p_c}}\geq [N(t)]^{\gamma}
\ee
for some universal constant $\gamma=\gamma(N,p)>0$ and for $t$ close enough to $T$. Indeed, we then conclude from (\ref{chohoeih}), (\ref{defvrescale}) and (\ref{defnt}) that:
$$|u(t)|_{L^{p_c}}=|v(0)|_{L^{p_c}}\geq [N(t)]^{\gamma}=|\log\lambda_u(t)|^{\gamma}\geq C|\log(T-t)|^{\gamma}$$ and (\ref{lowerboundcritical}) follows.\\

{\bf step 1} Propositions \ref{lemmarhonorm} and \ref{lemmarhonormd} apply to the rescaled solution.\\

We claim that for $t$ close enough to $T$:
\be
\label{estpropcleun}
\forall \tau_0\in [0,e^{N(t)}=\frac{1}{\lambda_u(t)}],  \ \ \mbox{Propositions \ref{lemmarhonorm} and \ref{lemmarhonormd} apply to} \ \  v(\tau_0).
\ee\\
We need to check that (\ref{estuimatestzero}), (\ref{estsizemzero}), (\ref{esttautildestar}) and (\ref{hypsupplemenatire}) hold.\\ 
Let $\tau_*=\frac{2}{\lambda_u(t)}\leq \tau_{max}$ from (\ref{defnt}) and $\tau_0\in [0,\frac{1}{\lambda_u(t)}]=[0,\frac{\tau_*}{2}]$ so that (\ref{esttautildestar}) holds. We have from (\ref{deflambdaut}):
\be
\label{estenergyev}
|\nabla v(0)|_{L^2}=1 \ \ \mbox{and} \ \ E(v_0)=\lambda_u(t)^{2(1-s_c)}E(u_0).
\ee 
In particular, $$\tau_*^{1-s_c}|E(v_0)|= \left(\tau_*\lambda^2_u(t)\right)^{1-s_c}|E(u_0)|= (2\lambda_u(t))^{1-s_c}|E(u_0)|\to 0 \ \ \mbox{as} \ \ t\to T$$ and (\ref{estuimatestzero}) follows. Note that (\ref{estenergyev}) implies $E(v_0)\to 0$ as $t\to T$ and thus $|\nabla v(0)|_{L^2}=1$ and (\ref{estggalibdeux}) imply (\ref{estsizemzero}) for $t$ close enough to $T$. Last, $$\lambda_v(\tau_0)=\left(\frac{1}{|\nabla v(\tau_0)|_{L^2}}\right)^{\frac{1}{1-s_c}}=\frac{\lambda_u(t-\lambda_u^2(t)\tau_0)}{\lambda_u(t)}$$ and thus from (\ref{estenergyev}): $$\lambda^{2(1-s_c)}_v(\tau_0)|E(v_0)|=\left(\lambda_v(\tau_0)\lambda_u(t)\right)^{2(1-s_c)}|E(u_0)|=\lambda_u(t-\lambda^2(t)\tau_0)|E(u_0)|\to 0 \ \ \mbox{as} \ \ t\to T$$ from $t\to T$ and $\lambda_u^2(t)\tau_0\leq \lambda_u(t)\to 0$ as $t\to T$. This concludes the proof of (\ref{estpropcleun}).\\

{\bf step 2} Construction of the channels.\\

Let $\alpha_2$ be the universal constant in (\ref{dispeestimategradient}). Let:
\be
\label{defmt}
M(t)=\frac{4|v(0)|_{L^{p_c}}}{C_{GN}}\geq 2, \ \ L(t)=\left[100[M(t)]^{\alpha_2}\right]^{\frac{1}{2(1-s_c)}}.
\ee
Observe that if $L(t)\geq e^{\frac{\sqrt{N(t)}}{2}}$, then (\ref{tobeprovedcrucial}) is proved so that we assume: 
\be
\label{estletdeux}
L(t)< e^{\frac{\sqrt{N(t)}}{2}}.
\ee
Following (\ref{defdzerofun}), let 
\be
\label{deffetey}
F(\tau)=\frac{\sqrt{\tau}}{\lambda_v(\tau)}.
\ee
We claim the existence of a family of channels: $\forall i\in [\sqrt{N(t)}, N(t)]$ integer, there exists $\tau_i\in [0,e^i]$ with 
\be
\label{condturctiondechaneesl}
F(\tau_i)\leq L(t)\ \ \mbox{and} \ \ \frac{1}{10L(t)}e^{\frac{i-1}{2}}\leq \lambda_v(\tau_i)\leq \frac{10}{L(t)}e^{\frac{i}{2}}.
\ee
Proof of (\ref{condturctiondechaneesl}): Pick $i\in [\sqrt{N(t)}, N(t)]$ and consider the set $${\cal{A}}=\{\tau\in[0,e^i] \ \ \mbox{with} \ \ \lambda_v(\tau)>\frac{e^{\frac{i+1}{2}}}{L(t)}\}.$$ Observe from (\ref{estenergyev}), (\ref{estletdeux}) that $$ \lambda_v(0)=1<\frac{e^{\frac{\sqrt{N(t)}}{2}}}{L(t)}\leq \frac{e^{\frac{i+1}{2}}}{L(t)}$$ and thus $0\notin{\cal{A}}$.\\
If ${\cal{A}}$ is non empty, the continuity of $\lambda_v(\tau)$ implies the existence of $\tau_i\in [0,e^i]$ such that $\lambda_v(\tau_i)=\frac{e^{\frac{i+1}{2}}}{L(t)}$ and then: $$F(\tau_i)=\frac{\sqrt{\tau_i}}{\lambda_v(\tau_i)}\leq L(t)\frac{\sqrt{e^i}}{e^{\frac{i+1}{2}}}\leq L(t)$$ and (\ref{condturctiondechaneesl}) is proved.\\
If ${\cal{A}}$ is empty, then equivalently: 
\be
\label{hoehofheo}
\forall \tau \in[0,e^i], \ \ \lambda_v(\tau)\leq \frac{e^{\frac{i+1}{2}}}{L(t)}\leq \frac{10}{L(t)}e^{\frac{i}{2}}.
\ee 
We now claim: 
\be
\label{fjp}
\exists\tau_i\in [e^{i-1},e^i] \ \ \mbox{with} \ \ F(\tau_i)\leq L(t).
\ee
Assume (\ref{fjp}), then $$\lambda_v(\tau_i)=\frac{\sqrt{\tau_i}}{F(\tau_i)}\geq \frac{\sqrt{e^{i-1}}}{L(t)}\geq \frac{e^{\frac{i-1}{2}}}{10L(t)}$$ and (\ref{condturctiondechaneesl}) follows from (\ref{hoehofheo}) and (\ref{fjp}).\\
It remains to prove (\ref{fjp}) which is a consequence of Proposition \ref{lemmarhonorm}. Indeed, if not, then: $$\forall \sigma\in [e^{i-1},e^i], \ \ F(\sigma)\geq L(t)\ \ \mbox{ie} \ \ |\nabla u(\sigma)|_{L^2}\geq \left(\frac{L(t)}{\sqrt{\sigma}}\right)^{1-s_c}$$ from (\ref{deflambdaut}) and (\ref{deffetey}). From (\ref{dispeestimategradient}), this implies:
\bee
& & [M(t)]^{\alpha_2}(e^i)^{1+s_c}  \geq   \int_0^{e^i}(e^i-\sigma)|\nabla v(\sigma)|_{L^2}^2d\sigma
\geq \int_{e^{i-1}}^{\frac{e^i}{2}}\sigma|\nabla v(\sigma)|_{L^2}^2d\sigma\\
& \geq & [L(t)]^{2(1-s_c)}\int_{e^{i-1}}^{\frac{e^i}{2}}\frac{\sigma}{\sigma^{1-s_c}}d\sigma> \frac{1}{100}(e^i)^{1+s_c}[L(t)]^{2(1-s_c)}>[M(t)]^{\alpha_2}(e^i)^{1+s_c}
\eee
from the definition (\ref{defmt}) of $L(t)$, and a contradiction follows which concludes the proof of (\ref{fjp}).\\

{\bf step 3} Uniform lower bound of the weighted $L^2$ norm on an annulus.\\

We now come to the heart of the proof. Let $i\in [\sqrt{N(t)}, N(t)]$ integer and $\tau_i$ be the times of the channels $\lambda_v(\tau_i)$ constructed from (\ref{condturctiondechaneesl}). We claim that there exist universal constants $\alpha_4(N,p), c_4(N,p)>0$ such that the following holds true: let the annulus 
\be
\label{deci}
{\cal{C}}_i=\{x, \ \ \frac{\lambda_v(\tau_i)}{[M(t)]^{\alpha_4}}\leq |x|\leq [M(t)]^{\alpha_4}\lambda_v(\tau_i)\},
\ee then:
\be
\label{weightedtgedibasic}
\forall i\in [\sqrt{N(t)}, N(t)], \ \ \int_{{\cal{C}}_i}|v(0)|^{p_c}\geq \frac{c_4}{[M(t)]^{\alpha_4s_cp_c}}.
\ee
Proof of (\ref{weightedtgedibasic}): From (\ref{estpropcleun}), we may apply (\ref{estsclaingldeux}) of Proposition \ref{lemmarhonormd} to $v(\tau_i)$. Now observe from (\ref{defmt}) and (\ref{condturctiondechaneesl}) that $F(\tau_i)\leq L(t)=[100[M(t)]^{\alpha_2}]^{\frac{1}{2(1-s_c)}}$ and thus the choice of $D_*$ given by (\ref{defdzero}) is uniform with respect to $i$ and satisfies $D^*(t)\leq [M(t)]^{\alpha_4}$ for some $\alpha_4(N,p)>0$. In other words, (\ref{estsclaingldeux}) implies:
\be
\label{estmowerbound}
\forall i\in[\sqrt{N(t)}, N(t)], \ \ \frac{1}{\lambda^{2s_c}_v(\tau_i)}\int_{|x|\leq [M(t)]^{\alpha_4}\lambda_v(\tau_i)}|v(0)|^2\geq c_3.
\ee
We now observe from H\"older inequality:
\bea
\label{cocoeiozohf}
\nonumber  & & \frac{1}{\lambda^{2s_c}_v(\tau_i)}\int_{|x|\leq \frac{\lambda_v(\tau_i)}{[M(t)]^{\alpha_4}}}|v(0)|^2  \leq \frac{1}{\lambda^{2s_c}_v(\tau_i)} \left(\int_{|x|\leq \frac{\lambda_v(\tau_i)}{[M(t)]^{\alpha_4}}}|v(0)|^{p_c}\right)^{\frac{2}{p_c}}\left(\frac{\lambda_v(\tau_i)}{[M(t)]^{\alpha_4}}\right)^{2s_c}\\
& \leq & \frac{CM^2(t)}{[M(t)]^{2s_c\alpha_4}}<\frac{c_3}{2}
\eea
for $\alpha_4>0$ large enough from (\ref{defmt}). We thus have from (\ref{estmowerbound}) and (\ref{cocoeiozohf}): 
\be
\label{nkek}
\forall i\in[\sqrt{N(t)}, N(t)], \ \ \frac{1}{\lambda^{2s_c}_v(\tau_i)}\int_{\frac{\lambda_v(\tau_i)}{[M(t)]^{\alpha_4}}\leq |x|\leq [M(t)]^{\alpha_4}\lambda_v(\tau_i)}|v(0)|^2\geq \frac{c_3}{2}.
\ee 
We now use H\"older again to estimate: 
$$\frac{1}{\lambda^{2s_c}_v(\tau_i)}\int_{\frac{\lambda_v(\tau_i)}{[M(t)]^{\alpha_4}}\leq |x|\leq C_4[M(t)]^{\alpha_4}\lambda_v(\tau_i)}|v(0)|^2\leq \frac{1}{\lambda^{2s_c}_v(\tau_i)}\left(\int_{{\cal{C}}_i}|v(0)|^{p_c}\right)^{\frac{2}{p_c}}\left(\lambda_v(\tau_i)[M(t)]^{\alpha_4}\right)^{2s_c}$$ which together with (\ref{nkek}) concludes the proof of (\ref{weightedtgedibasic}).\\

{\bf step 4} Conclusion.\\

We are now in position to conclude the proof of Theorem \ref{thmmain}. Let $p(t)>1$ be an integer such that 
\be
\label{dep} 
1000[M(t)]^{2\alpha_4}\leq e^{\frac{p(t)}{2}}\leq 4000^2[M(t)]^{2\alpha_4}.
\ee
We consider two cases.\\ 
If $p(t)>\sqrt{N(t)}$ then $|v(0)|_{L^{p_c}}^{2\alpha_4}=C[M(t)]^{2\alpha_4}\geq Ce^{\frac{p(t)}{2}}\geq Cp(t)\geq C\sqrt{N(t)}$ and (\ref{tobeprovedcrucial}) is proved.\\ 
If
\be
\label{dhoezhooez}
p(t)<\sqrt{N(t)},
\ee
we observe from the decoupling (\ref{condturctiondechaneesl}) and the choice (\ref{dep}) that: $$\frac{\lambda_v(\tau_{i+p})}{[M(t)]^{\alpha_4}}\geq \frac{1}{[M(t)]^{\alpha_4}}\frac{e^{\frac{i+p-1}{2}}}{10L(t)}> [M(t)]^{\alpha_4}\frac{10e^{\frac{i}{2}}}{L(t)}\geq [M(t)]^{\alpha_4}\lambda_v(\tau_i)$$ and thus the family of annuli ${\cal{C}}_i$ given by (\ref{deci}) satisfies: $$\forall (i,i+p)\in[\sqrt{N(t)}, N(t)]^2, \ \ {\cal{C}}_{i+p}\cap{\cal{C}}_{i}=\emptyset.$$ From (\ref{dhoezhooez}), we conclude that there are at least $\frac{N(t)}{10p(t)}\geq \frac{1}{10}\sqrt{N(t)}$ disjoint annuli on which the uniform lower bound (\ref{weightedtgedibasic}) holds. Summing over these annuli yields:
$$
|v(0)|^{p_c}_{L^{p_c}}\geq \Sigma_{k=0}^{\frac{N(t)}{10p(t)}}\int_{{\cal{C}}_{\sqrt{N(t)}+kp(t)}} |v(0)|^{p_c}\geq C\frac{\sqrt{N(t)}}{[M(t)]^{\alpha_4p_cs_c}}=C\frac{\sqrt{N(t)}}{|v(0)|_{L^{p_c}}^{\alpha_4p_cs_c}}$$ which concludes the proof of (\ref{tobeprovedcrucial}).\\

This concludes the proof of Theorem \ref{thmmain}.


\section*{Appendix A: Lifetimes in $\dot{H}^{s_c}$ and $\dot{H}^1\cap\dot{H}^{s_c}$}


Let $s_c<s\leq 1$, $u_0\in \dot{H}^{s_c}\cap \dot{H}^s$ and $u(t)$ the corresponding solution to (\ref{nls}) with $T_{\dot{H}^{s_c}\cap \dot{H}^s}$, $T_{\dot{H}^{s_c}}$ its life times given by the local Cauchy theory in respectively $\dot{H}^{s_c}\cap \dot{H}^s$ and $\dot{H}^{s_c}$, then standard arguments based on Strichartz estimates ensure that:
\be
\label{estaproeubeho}
T_{\dot{H}^{s_c}\cap \dot{H}^s}=T_{\dot{H}^{s_c}}.
\ee 
For the sake of completeness, we present a simple argument for $s=1$ based on the conservation of the energy.\\
By definition, $T_{\dot{H}^{s_c}\cap \dot{H}^1}\leq T_{\dot{H}^{s_c}}$. Let now $u_0\in \dot{H}^{s_c}\cap \dot{H}^1$ so that $T_{\dot{H}^{s_c}\cap \dot{H}^1}<+\infty$ or equivalently 
\be
\label{ceoehoz}
\lim_{t\to T_{\dot{H}^{s_c}\cap \dot{H}^1}}|\nabla u(t)|_{L^2}= +\infty.
\ee
 We claim that on any subsequence $t_n\to T_{\dot{H}^{s_c}\cap \dot{H}^1}$, $u(t_n)$ does not have a strong limit in $\dot{H}^{s_c}$, what thus implies (\ref{estaproeubeho}) from the local Cauchy theory in $\dot{H}^{s_c}$.\\ 
By contradiction, if $u(t_n)$ is compact in $\dot{H}^{s_c}$, then in particular no concentration occurs in $\dot{H}^{s_c}$ so that given $\e>0$ small enough, we may find a partition of $\R^N$ with $M(\e)$ balls $B_i=B(x_i,r(\e))$ and corresponding localizing cut-off $\chi_i=\chi(\frac{x-x_i}{r(\e)})$  such that:
$$
\forall n\geq N(\e), \ \ \forall i\in [0,M(\e)], \ \ |\chi_i u(t_n)|_{\dot{H}^{s_c}}<\e.
$$
We then may localized the Gagliardo-Nirenberg inequality (\ref{galiehozinc}) as follows -where $C$ is a large constant depending on $N,p$-:
\bee
\int |u(t_n)|^{p+1} & \leq  & C\Sigma_{i=0}^{M(\e)}|\chi_i u(t_n)|_{L^{p+1}}^{p+1}\leq C\Sigma_{i=0}^{M(\e)}|\chi_i u(t_n)|_{\dot{H}^{s_c}}^{p-1}|\chi_iu(t_n)|_{\dot{H}^1}^2\\
& \leq  & C\e^{p-1}|\nabla u(t_n)|^2_{L^2}+C\e^{p-1}\Sigma_{i=0}^{M(\e)}|\nabla \chi_i u(t_n)|_{L^2}^2\\
& \leq  & C\e^{p-1}|\nabla u(t_n)|^2_{L^2}+C_{\e}\Sigma_{i=0}^{M(\e)}|u(t_n)|_{L^{p_c}(B_i)}^2\\
& \leq &  C\e^{p-1}|\nabla u(t_n)|^2_{L^2}+C_{\e}|u(t_n)|_{L^{p_c}}^{p_c}+C_{\e}\\
& \leq  & C\e^{p-1}|\nabla u|^2_{L^2}+C_{\e}|u(t_n)|^{p_c}_{\dot{H}^{s_c}}+C_{\e}
\eee
where we used H\"older inequality with $p_c=\frac{N(p-1)}{2}=\frac{2N}{N-2s_c}>2$ and the Sobolev embedding $\dot{H}^{s_c}\subset L^{p_c}$. From the conservation of the energy, this implies for $\e$ small enough that $$\limsup_{t_n \to T_{\dot{H}^{s_c}\cap \dot{H}^1}}|\nabla u(t_n)|_{L^2}<C_{\e}<+\infty$$ which contradicts (\ref{ceoehoz}). See \cite{MR6} for similar arguments. This concludes the proof of (\ref{estaproeubeho}).

\end{document}